\documentclass[a4paper,11pt]{amsart}
\DeclareRobustCommand{\SkipTocEntry}[5]{}

  \usepackage{ tipa }

  \usepackage{rotating}
  \usepackage{floatpag}
  \rotfloatpagestyle{empty}

 \usepackage{amsmath}% if you are using this package,
  \usepackage{amsthm}
  \usepackage{graphicx}

  \usepackage{multind}
  \makeindex{authors}
  \makeindex{subject}

  \theoremstyle{plain}% default
  \newtheorem{theorem}{Theorem}[section]
  
  \newtheorem{proposition}[theorem]{Proposition}

  \theoremstyle{definition}
  \newtheorem{definition}[theorem]{Definition}
  \newtheorem{example}[theorem]{Example}
   \newtheorem{question}[theorem]{Question}

  \theoremstyle{remark}
  \newtheorem*{remark}{Remark}

\usepackage{amsmath }
\usepackage{marvosym}

\usepackage[T1]{fontenc}
\usepackage{empheq}
\usepackage{relsize}
\usepackage{amsmath}
\usepackage{bm}

\usepackage{upgreek}
\usepackage{ esint }
\usepackage{color}
\usepackage{comment}
\usepackage{amssymb}
\usepackage{amsfonts}
\usepackage{graphicx}

\usepackage{amscd}

\usepackage{slashed}
\usepackage{pdflscape}
\usepackage{enumerate}
\usepackage{multirow}
\usepackage{stmaryrd}

\usepackage{tikz-cd}
\usepackage{tikz}
\usepackage{tkz-graph}
\usetikzlibrary{topaths}
\usetikzlibrary{arrows}
\usepackage{mathdots}
\usetikzlibrary{shapes,snakes}
\usetikzlibrary{decorations.shapes}
\usetikzlibrary{shapes,decorations,shadows}
\usetikzlibrary{decorations.pathmorphing}
\usetikzlibrary{decorations.shapes}
\usetikzlibrary{fadings}
\usetikzlibrary{patterns}
\usetikzlibrary{calc}
\usetikzlibrary{decorations.text}
\usetikzlibrary{decorations.footprints}
\usetikzlibrary{decorations.fractals}
\usetikzlibrary{shapes.gates.logic.IEC}
\usetikzlibrary{shapes.gates.logic.US}
\usetikzlibrary{fit,chains}
\usetikzlibrary{positioning}
\usepgflibrary{shapes}
\usetikzlibrary{scopes}
\usepackage{bbold}
\usepackage[linkcolor=blue]{hyperref}
\usepackage{mathrsfs}
\usepackage[colorinlistoftodos]{todonotes}

\usepackage[T1]{fontenc}

\definecolor{blue}{rgb}{.255,.41,.884} % RoyalBlue of svgnames
\definecolor{red}{rgb}{1, 0, 0} % Red of svgnames
\definecolor{green}{rgb}{.196,.804,.196} % LimeGreen of svgnames
\definecolor{yellow}{rgb}{1,.648,0} % Orange of svgnames
\definecolor{pink}{rgb}{1,0.5,0.5}

\usepackage{pdfsync}

\usepackage{graphicx} 
\usepackage{amsmath} 
\usepackage{amsfonts}
\usepackage{amssymb}

\newcommand{\be}{\begin{equation}}

\newcommand{\ee}{\end{equation}}

\newcommand{\Odane}{\O}

\newcommand{\om}{\omega}

\newcommand{\cB}{{\mathcal B}}
\newcommand{\cN}{\mbox{$\mathcal{N}$}}

\newcommand{\II}{{ I\hspace{-.8mm}I}}
\newcommand{\IIo}{\mathring{\!{ I\hspace{-.8mm} I}}{\hspace{.2mm}}}

\newcommand{\si}{\sigma}

\newcommand{\ba}{\begin{array}}

\newcommand{\ea}{\end{array}}

\newcommand{\beq}{\begin{eqnarray}}

\newcommand{\eeq}{\end{eqnarray}}

\newtheorem{lm}{lemma}

\newtheorem{thee}{theorem}

\newtheorem{proo}{proposition}

\newtheorem{co}{corollary}

\newtheorem{rem}{remark}

\newtheorem{deff}{definition}

\newcommand{\bd}{\begin{deff}}

\newcommand{\ed}{\end{deff}}

\newcommand{\bl}{\begin{lm}}

\newcommand{\el}{\end{lm}}

\newcommand{\bp}{\begin{proo}}

\newcommand{\ep}{\end{proo}}

\newcommand{\bt}{\begin{thee}}

\newcommand{\et}{\end{thee}}

\newcommand{\bc}{\begin{co}}

\newcommand{\ec}{\end{co}}

\newcommand{\brm}{\begin{rem}}

\newcommand{\erm}{\end{rem}}

\usepackage{t1enc}

\def\Cal{\mathcal}

\newcommand{\bS}{\mathbb{S}}

\newcommand{\newc}{\newcommand}

\let\ccdot.
\def\cmdot{\hbox to 2.5pt{\hss$\hh\cdot\hh$\hss}}

\newc{\aR}{\mbox{\boldmath{$ R$}}}
\newc{\aS}{\mbox{\boldmath{$ S$}}}
\newc{\aT}{\mbox{\boldmath{$ T$}}}
\newc{\aW}{\mbox{\boldmath{$ W$}}}

\newc{\aD}{\mbox{\boldmath{$ D$}}\hspace{-.2mm}}

\newc{\aK}{\mbox{\boldmath{$ K$}}}
\newc{\aL}{\mbox{\boldmath{$ L$}}}
\newcommand{\ce}{{\Cal E}}

\newcommand{\ct}{{\Cal T}}

\newcommand{\bT}{{\Bbb T}}

\usepackage{amssymb}
\usepackage{amscd}

\newcommand{\nn}[1]{(\ref{#1})}

\newcommand{\bg}{\mbox{\boldmath{$ g$}}}

\let\t=\tau

\newc{\obstrn}[2]{B^{#1}_{#2}}

\newcommand{\rpl}                         % +) or <+
{\mbox{$
\begin{picture}(12.7,8)(-.5,-1)
\put(0,0.2){$+$}
\put(4.2,2.8){\oval(8,8)[r]}
\end{picture}$}}

\newcommand{\lpl}                         % (+ or +>
{\mbox{$
\begin{picture}(12.7,8)(-.5,-1)
\put(2,0.2){$+$}
\put(6,2.7){\oval(7,7)[l]}
\end{picture}$}}

\usepackage{ifthen}

\newc{\tensor}[1]{#1}
\newc{\Mvariable}[1]{\mbox{#1}}
\newc{\down}[1]{{}_{#1}}
\newc{\up}[1]{{}^{#1}}

\newc{\JulyStrut}{\rule{0mm}{6mm}}
\newc{\midtenPan}{\mbox{\sf S}}
\newc{\midten}{\mbox{\sf T}}
\newc{\midtenEi}{\mbox{\sf U}}
\newc{\ATen}{\mbox{\sf E}}
\newc{\BTen}{\mbox{\sf F}}
\newc{\CTen}{\mbox{\sf G}}

\def\sideremark#1{\ifvmode\leavevmode\fi\vadjust{\vbox to0pt{\vss
 \hbox to 0pt{\hskip\hsize\hskip1em
 \vbox{\hsize2cm\tiny\raggedright\pretolerance10000
  \noindent #1\hfill}\hss}\vbox to8pt{\vfil}\vss}}}

\numberwithin{equation}{section}

\newcommand{\hh}{{\hspace{.3mm}}}

\newcommand{\cc}{\boldsymbol{c}}

\newcommand{\sss}{\scriptscriptstyle}

\newcommand{\pdot}{{\textstyle\boldsymbol \cdot}\hspace{.05mm}}

\renewcommand\geq{\geqslant}
\renewcommand\leq{\leqslant}

\newcommand{\ext}{{\rm d}}

\begin{document}

\vspace{3pt}

\renewcommand{\arraystretch}{1}

\renewcommand{\today}{}
\title{Singular Yamabe and Obata Problems
}
\author{ A. Rod Gover${}^\sharp$ \&  Andrew Waldron${}^\natural$}

\address{${}^\sharp$
  Department of Mathematics\\
  The University of Auckland\\
  Private Bag 92019\\
  Auckland 1142\\
  New Zealand,  and\\
  Mathematical Sciences Institute, Australian National University, ACT 
  0200, Australia} \email{gover@math.auckland.ac.nz}
  
  \address{${}^{\natural}$
  Center for Quantum Mathematics and Physics (QMAP)\\
  Department of Mathematics\\ 
  University of California\\
  Davis, CA95616, USA} \email{wally@math.ucdavis.edu}

\vspace{3pt}

\renewcommand{\arraystretch}{1}

\begin{abstract}

  A conformal geometry determines a distinguished,
    potentially singular, variant of the usual Yamabe problem, where the
    conformal factor can change sign. When a smooth 
    solution does change sign, its zero locus is a smoothly
    embedded separating hypersurface that, in dimension three, is necessarily a Willmore
    energy minimiser or, in higher dimensions, satisfies a  conformally invariant analog of the Willmore
    equation.  In any case the zero locus is critical for a conformal
    functional that generalises the total $Q$-curvature by including
    extrinsic data.  These observations lead to some interesting
    global problems that include  natural singular
    variants of a classical problem solved by Obata.

\end{abstract}

\vspace{10cm}

%\begin{abstract}
%
%Stuff

%\vspace{10cm}
%
%\noindent
%%\begin{keywords}
%{\sf \tiny Keywords: 
%Conformal Geometry, Embedded Manifolds with Boundary, 
%Extrinsic Laplacian Powers and Boundary Operators,
%$Q$ and~$T$-transgression Curvatures,
%Minimal Hypersurface Asymptotics,
% AdS/CFT,  Anomalies,  Renormalized Volume, 
%     Conformally Compact, Yamabe problem, Willmore Energies with Boundary \color{black}}
% %\end{keywords}

%\end{abstract}

%58E30  	Variational principles
%49S05  	Variational principles of physics
%51P05  	Geometry and physics
%53A30  	Conformal differential geometry
%53A55  	Differential invariants (local theory), geometric objects
%53C80  	Applications to physics (Global differential geoemtry)
%53B50  	Applications to physics (Local differential geometry)
%53A07  	Higher-dimensional and -codimensional surfaces in Euclidean~$n$-space
%53C21  	Methods of Riemannian geometry, including PDE methods; curvature restrictions
%53B15  	Other connections
%83C99  	None of the above, but in this section (General relativity)

\maketitle

\pagestyle{myheadings} \markboth{Gover \& Waldron}{Singular Obata Problem}

\newpage

%\tableofcontents

\newcommand{\balpha}{{\bm \alpha}}
\newcommand{\balphas}{{\scalebox{.76}{${\bm \alpha}$}}}
\newcommand{\bnu}{{\bm \nu}}
\newcommand{\bnus}{{\scalebox{.76}{${\bm \nu}$}}}
\newcommand{\bnuss}{\hh\hh\!{\scalebox{.56}{${\bm \nu}$}}}

\newcommand{\bmu}{{\bm \mu}}
\newcommand{\bmus}{{\scalebox{.76}{${\bm \mu}$}}}
\newcommand{\bmuss}{\hh\hh\!{\scalebox{.56}{${\bm \mu}$}}}

\newcommand{\btau}{{\bm \tau}}
\newcommand{\btaus}{{\scalebox{.76}{${\bm \tau}$}}}
\newcommand{\btauss}{\hh\hh\!{\scalebox{.56}{${\bm \tau}$}}}

\newcommand{\bsigma}{{\bm \sigma}}
\newcommand{\bsigmas}{{{\scalebox{.8}{${\bm \sigma}$}}}}
\newcommand{\bbeta}{{\bm \beta}}
\newcommand{\bbetas}{{\scalebox{.65}{${\bm \beta}$}}}

\renewcommand{\bS}{{\bm {\mathcal S}}}
\newcommand{\bB}{{\bm {\mathcal B}}}
\newcommand{\bC}{{\bm {\mathcal C}}}

\renewcommand{\bT}{{\bm {\mathcal T}}}
\newcommand{\bM}{{\bm {\mathcal H}}}

\newcommand{\go}{{\mathring{g}}}
\newcommand{\nuo}{{\mathring{\nu}}}
\newcommand{\alphao}{{\mathring{\alpha}}}

\newcommand{\Ell}{\mathscr{L}}
\newcommand{\density}[1]{[g\, ;\, #1]}

\newcommand{\ceedot}{{\scalebox{2}{$\cdot$}}}

\newcommand{\Langle}{{\bm \langle}}
\newcommand{\Rangle}{{\bm \rangle}}
\newcommand{\Lodz}{\widehat{\mathbf L}}
\newcommand{\Dhat}{\mbox{\bf \DJ}}

\section{Introduction}

On a closed Riemanian  manifold $(M,\bar{g})$ the Yamabe problem concerns
finding a conformally related metric $g=e^{2\omega}$, for some $\om\in C^\infty (M)$, that has constant
scalar curvature.  The statement and ultimate solution of this problem by Yamabe, Trudinger, Aubin,
and Schoen~\cite{Yamabe, Trudinger, Aubin,  Schoen1} were milestones in conformal geometry and
geometric analysis. Another---apparently unrelated---problem  of great interest has
been that of finding critical points of the Willmore energy~\cite{Willmore,Marques}. The Willmore energy $\mathcal W$ for a
closed surface~$\Sigma$ in Euclidean 3-space  is given by
$$
\mathcal{W}=\int_\Sigma (H^2-K)\, , 
$$
where $H$ is the surface mean curvature, and $K$ its Gauss
curvature. A key property of this energy is that it is invariant under
conformal transformations~\cite{WillmoreB}.
We shall call the Euler--Lagrange equation for this  energy the {\it Willmore equation}, this equation determines a conformally invariant quantity termed the {\it Willmore invariant}.

Recently it has been observed that there is a rather interesting link
between natural generalisations of these problems. For that 
one replaces the usual Yamabe problem with a singular variant---namely on a  Riemannian manifold $(M,\bar{g})$ of dimension $d$ one seeks a smooth (here and throughout smooth means $C^\infty$) function
$\si$ such that
\begin{equation}\label{sYam}
|\ext \sigma|_{\bar g}^2 - \frac{2\sigma }d\,  
 \Big(\Delta^{\bar g} \sigma +\frac\sigma{2(d-1)} \,  \, Sc^{\bar g}\Big) =1\,  .
\end{equation}
This is the usual Yamabe-type equation requiring that the scalar curvature
satisfies $\operatorname{Sc}^{g}=-d(d-1)$, where $g=\si^{-2}\bar{g}$,
except that we allow the possibility that $\si$ changes sign. In the
latter case the metric $g$ is singular along the zero locus~$\mathcal{Z}(\si)$ of $\si$, and it is evident from the equation and
setup here, that this zero locus is an embedded hypersurface
({\it i.e.}, codimension-1 submanifold) with an induced conformal
structure. See Proposition~\ref{sep-prop} below for more detail. More striking is that  equation~\nn{sYam} puts an interesting
restriction on the conformal embedding of $\mathcal{Z}(\si)$; if $d=3$
it is necessarily Willmore ({\it i.e.}, a Willmore energy minimiser) and in
higher dimensions this condition {\em defines} a conformally invariant
generalisation of the  Willmore equation; see Theorem
\ref{mainonex} which follows~\cite{Andersson,GW13,GW15}. It turns out that
this higher Willmore equation is the Euler-Lagrange
equation (with respect to variation of embedding) of an action~\cite{Grahamnew}, which in fact can be expressed as an integral of a
quantity that generalises, by the addition of extrinsic curvature
terms~\cite{GW14,GW16a,GAW}, the Branson $Q$-curvature of \cite{BQ,BO} (see the review~\cite{WhatQ}). 

The picture just described captures some of the important local
aspects of the link between (higher) Willmore minimisers and solutions
to the singular Yamabe equation~\nn{sYam}. A key purpose of the
current note is to  point out some very interesting features and questions linked to the 
global version of this problem on closed manifolds. This culminates in the main
questions which we introduce in Section~\ref{Qsect}. In Section
\ref{background} we review briefly some of the background and mention
informally a singular variant of one of the Obata problems.  Section
\ref{tractors} describes some key tools from conformal tractor
calculus required to handle hypersurfaces embedded in conformal manifolds.

\section{Background and a singular Obata problem} \label{background}

Let $(M^d,\cc)$ be a closed, orientable, Riemannian signature, conformal manifold where $\cc\ni[g]=[\Omega^2 g]$ denotes an equivalence class of smooth, conformally related metrics with 
$0<\Omega\in C^\infty M$, and $d\geq 3$. The trace-free Schouten tensor of $g\in \cc$ defined by
$$\mathring P^g =\frac1{d-2} \, \big(Ric^g -\frac 1d\,  g\,  Sc^g\big)\, , $$ 
obeys the conformal transformation law
$$
\mathring P_{ab}^{\Omega^2 g} = 
\mathring P^g_{ab} 
- \nabla^g{}_{(a} \Upsilon_{b)\circ}
 + \Upsilon_{(a}\Upsilon_{b)\circ}\, .
$$
Here $\nabla$ is the Levi-Civita connection, $\Upsilon:=\ext \log \Omega$, and we have employed an abstract index notation to denote sections of tensor bundles over $M$ as well as the notation $X_{(ab)\circ}:=X_{ab}-\frac1d \, g_{ab} \langle g, X\rangle_g$
for projection to the trace-free part of $X\in \Gamma(\odot^2 T^*M)$; this projection is  independent of the choice of $g\in \cc$. Also note that
$\langle X,Y\rangle_g:= g^{ac}g^{bd}
X_{ab}
Y_{cd}$ and $|X|_g^2:=\langle X,X\rangle$.
Metrics whose trace-free Schouten 
tensor vanishes are called {\it Einstein}.
\medskip

Now, supposing that $g$ has constant scalar curvature, the Bianchi identity gives (suppressing the $g$ dependence)
$$
\nabla^a \mathring P_{ab}=0\, .
$$
Thus, in that case it follows that
$$
\int_{(M,g)} \mathring P^{ab} \nabla_a \nabla_b \Omega^{-1} = 0 \Rightarrow
\int_{(M,g)} \Omega^{-1}\mathring P^{ab} \big(\nabla_a \Upsilon_b - \Upsilon_a\Upsilon_b\big)=0\, . 
$$
Calling $\Omega^2 g=\bar g$, we see that
$ 
\int_{(M,g)}\Omega^{-1} \langle \mathring P^g,\mathring P^{\bar g}\rangle=
\int_{(M,g)}\Omega^{-1} \, |\mathring P^g|^2
$.
%Noting that  the Schwarz inequality implies
%$
%\big(\int_{(M,g)} \langle \mathring P^g,
%\mathring P^{\bar g}\rangle\big)^2$ $\leq$ 
%\int_{(M,g)} |\mathring P^g|^2 \, .\, 
%\int_{(M,g)} |\mathring P^{\bar g}|^2 
%$,
The Schwarz inequality then implies that
% \edz{THink about factors $\Omega$} here we have
$$
\int_{(M,g)}\Omega^{-1} |\mathring P^g|_g
^2\leq
\int_{(M,g)}\Omega^{-1} |\mathring P^{\bar g}|_g^2\, .
$$
The above simple argument due to Schoen~\cite{Schoen}, implies that if two metrics are conformally related and one is Einstein, while the other has constant scalar curvature, then they must in fact both be Einstein. This proves the following classical theorem due to Obata~\cite{Obata}.

\begin{theorem}\label{OBATA}
Let $(S^d,\cc_{\rm round})$ be the sphere equipped with its standard,  conformally flat, class of metrics. Then
if $g\in \cc_{\rm round}$ has constant scalar curvature, $g$ must be Einstein. 
\end{theorem}

It is easy to construct an example of this phenomenon. 
\begin{example}
Let $h:S^d \to [-1,1]$ be the height function on the standard sphere in~${\mathbb R}^{d+1}$,
$$
S^d:=\{x^2+y^2+\cdots + z^2=1\}\, ,
$$
given by $h(x,y,\ldots,z)=z$. Then
if $\bar g$ is the standard sphere metric, 
an elementary computation shows that, where it is defined,  the conformally related metric 
$$
g=\frac{\bar g}{(h-k)^2}\, ,
$$
 has constant scalar curvature
$$
Sc^g=d(d-1)(k^2-1)\, .
$$
%The Weyl tensor of $g$ necessarily vanishes. 
We give a simple argument in Section~\ref{sphere} that, for $k>1$, the
 metric $g$ is indeed Einstein and, since its Weyl tensor necessarily vanishes,
 isometric to the standard sphere metric.
\end{example}

%COMMENT USE SCHOEN, EXPLICT r DIFFEO OR stereo to EXPLAIN BELOW

%\noindent
An intriguing feature of this example is that 
when $k=1$ the metric~$g$ becomes singular at the north pole $z=1$. However, 
in that case, the scalar curvature vanishes and $(S^d/\{z=1\},g)$ is isometric to Euclidean space.
 This can be checked explicitly by writing $\bar g = d\theta^2 + \sin^2\theta \, ds^2_{S^{d-1}}$, where $h=z=\cos\theta$, and then noting that the change of coordinates $r=\sin\theta/(1-\cos\theta)$ gives the required isometry between $g$ and the flat metric $dr^2 + r^2 ds^2_{S^{d-1}}$. Geometrically this is the stereographic projection of $S^d$ to the hyperplane $z=1$.

%\edz{RQ: Is the isometry\\
% to the flat metric\\
%  already evident?\\
%   What is the logic here?\\
%YES, either just compute\\
% the curvatures or\\ 
% note that the change\\
%  of variable
%$
%r=1/(1-\cos(theta)
%$\\
%brings the 
%metric to\\
%$g=dr^2 + r^2 ds^2_{S^{d-1}}$\\
%there must be\\
%a similar triang.s\\
%stereo-proj picture...
%}

When $|k|<1$, $g$ becomes singular along a hypersurface of height $k$.
%Focussing (for simplicity) on the case $k=0$, 
The Riemannian manifold $(S_+,g_+)$, given by the above data restricted to  $z>k$, is isometric to the Poincar\'e ball.  
%\edz{RQ: How do we see ``is isometric to the Poincar\'e disc''?\\
%Perhaps we should leave as a claim that will be established in Section \ref{sphere} -- I have added a proof there. Sounds good}
The same applies to
%the southern hemisphere~
$(S_-,g_-)$ with~$z<k$. Observe that  the conformal rescaling function~$\Omega=\frac{1}{h-k}$ is
singular along the hypersurface $\Sigma=S^{d-1}$ at $z=k$, and changes sign  in~$S_-$ because the function $h-k$ is a defining function  for the hypersurface embedding $\Sigma\hookrightarrow S^d$.
These facts are recovered in
Section~\ref{sphere} from a different perspective that unifies the three cases discussed above and provides information about the geometry of the embedding $\Sigma\hookrightarrow M$.

\medskip
 Given a Riemannian manifold $(M^d,\bar g)$ and a smooth function~$\sigma$ we define
 \begin{equation}\label{scurvy}
 S(\bar g,\sigma) := |\ext \sigma|_{\bar g}^2 - \frac{2\sigma }d\,  
 \Big(\Delta^{\bar g} \sigma +\frac\sigma{2(d-1)} \,  \, Sc^{\bar g}\Big)\, ,
 \end{equation}
 which obeys 
$$
S(\Omega^2 \bar g,\Omega\sigma) = S(\bar g,\sigma)\, ,
$$
for any $0<\Omega\in C^\infty M$.
When the part of the jet of $\sigma$ given by $\big(\sigma,\ext  \sigma , -\tfrac1d (\Delta^{\bar g} \sigma +\tfrac{\sigma}{2(d-1)}Sc^{\bar g})\big)$ is nowhere vanishing, $S(\bar g, \sigma)$ has  a natural interpretation as a curvature that we term the {\it $S$-curvature}. Indeed, 
 away from the zero locus $\Sigma$ of $\sigma$,
\begin{equation}\label{Sc}
S(g,1)=-\frac{Sc^g}{d(d-1)} \mbox{ where } g=\sigma^{-2} \bar g\, ,
\end{equation}
so the $S$-curvature smoothly extends the scalar curvature of the singular metric 
$g=\sigma^{-2}\bar g$ 
%that is only defined on $M\backslash \Sigma$, 
to all of $M$.

Given its nature as a curvature that extends (and generalises) the
scalar curvature, it is interesting to consider the problem of finding
functions $\si$ that yield $S(\bar{g},\si)$ constant. If $\si$ has
fixed sign then this  boils down to the usual Yamabe problem on
closed manifolds. However, in general other solutions are possible and
hence it is interesting to consider the case that $\si$ has a zero
locus. Note that in view of (\ref{scurvy}), a non-trivial zero locus is
impossible if $S(\bar{g},\si)$ is strictly negative. If
$S(\bar{g},\si)=0$ then only isolated zeros are possible~\cite{Goal}. But
observe that when $S(\bar g, \sigma)$ is strictly positive, the zero
locus of $\sigma$, if non-empty, is a hypersurface with defining
function $\sigma$ (see~\cite{Goal}). Here we shall be most interested in this last case,
 and  if $\sigma$   satisfies
 \begin{equation}\label{singyam}
S(\bar g,\sigma)=1\, ,
\end{equation}
then the singular metric $g=\sigma^{-2} g$ has (by virtue of
Equation~\nn{Sc}) constant negative scalar curvature $-d(d-1)$ on
$M\backslash \Sigma$. Hence the above display is called the {\it
  singular Yamabe equation} (see~\cite{Mazzeo}).  For $\Sigma$ closed and orientable, it
is known that one-sided solutions for $\sigma$ exist~\cite{Loewner, Aviles, Mazzeo, Andersson}.

This brings us to a distinct but
closely related local problem. Namely given an embedded hypersurface
$\Sigma$, can we find a local defining function for~$\Sigma$ that is
distinguished in the sense that $\si$ solves $S(\bar{g},\si)=1$. It
turns out that for smooth solutions $\si$, this is obstructed by an
interesting conformal invariant of the extrinsic geometry.  For
example for surfaces~$\Sigma$ we have the following result, which
follows from the work of~\cite{Andersson} and an observation made
in~\cite{GW13,GW15}.

%% question whether it is possible to find $\sigma$ such that the
%% $S$-curvature is a positive constant and $\sigma$ is defining for a
%% given hypersurface is particularly interesting, especially in light of
%% the following result which follows from the work of~\cite{ACF} and an
%% observation made in~\cite{US}.
%%  % 
%  an embedded hypersurface $\Sigma\hookrightarrow M$ with defining function $\sigma\in C^\infty M$ (so $\Sigma = {\mathcal Z}(\sigma)$ and $\ext  \sigma|_\Sigma\neq 0$\edz{define in intro!!}), the {\it $S$-curvature}
% defined by
% \begin{equation}\label{scurvy}
% S(\bar g,\sigma) := |\ext \sigma|_{\bar g}^2 - \frac{2\sigma }d\,  
% \Big(\Delta^{\bar g} \sigma +\frac\sigma{2(d-1)} \,  \, Sc^{\bar g}\Big)\, ,
% \end{equation}
%obeys
%$$
%S(\Omega^2 \bar g,\Omega\sigma) = S(\bar g,\sigma)\, ,
%$$
%for any $0<\Omega\in C^\infty M$.
%Moreover, away from the zero locus $\Sigma$ of $\sigma$,
%$$
%S(g,1)=-\frac{Sc^g}{d(d-1)} \mbox{ where } g=\sigma^{-2} \bar g\, .
%$$
%Hence, the $S$-curvature smoothly extends the scalar curvature of the singular metric \edz{define singular metric in intro} 
%$g=\sigma^{-2}\bar g$ 
%%that is only defined on $M\backslash \Sigma$, 
%to all of $M$. 

\begin{proposition}\label{Willmore}
  Let $(M^3,\bar g)$ be a Riemannian three-manifold. Then a surface $\Sigma$ admits a defining function satisfying
  %\begin{equation}\label{nottoB}
$$
S(\bar g,\sigma)=1 + \sigma^{4} T\, ,
$$
%\end{equation}
for some $T\in C^\infty M$, if and only if
it is Willmore.
\end{proposition}

\noindent
The above proposition is a special case of a rather uniform picture, from~\cite{Andersson,GW13,GW15,GW161}:
\begin{theorem}\label{mainonex}
Given a hypersurface $\Sigma\hookrightarrow (M^d,\bar g)$ embedded in
a Riemannian $d$-manifold, there
exists a smooth defining function  such that
$$
S(\bar g, \sigma)=1 +  \sigma^d B\,   ,
$$
where the {\it obstruction}
$B_\Sigma=B|_\Sigma$ is a hypersurface conformal invariant and so, in particular, depends only on the
conformal embedding $\Sigma\hookrightarrow (M^d,[\bar g])$.
\end{theorem}

Conformal hypersurface invariants are defined in \cite{GW15}. The quantity~$B_\Sigma$ obstructs smooth solutions $\sigma$ to the equation $ S(\bar g,
\sigma)=1 $ subject to the boundary condition that $\Sigma$ is the
zero locus of $\sigma$.  When $d=3$, the obstruction $B_\Sigma$ is the
functional gradient of the Willmore energy functional, and thus defines the Willmore invariant discussed above. It has leading
term $-\frac{1}{3}\Delta_\Sigma H$, where $H$ is the hypersurface
mean curvature. For hypersurfaces of higher dimension the
obstruction~$B_\Sigma$ generalises this and, in particular, for $d$ even  has leading term
$$
\Delta^{\frac{d-1}{2}}_\Sigma H, \quad\mbox{up to a constant}.
$$
Thus we term $B_{\Sigma}$ a {\it higher Willmore invariant}.

It turns out that the invariants $B_\Sigma$ are variational: in $d$
dimensions, generalized Willmore energy functionals exist with
functional gradient $\cB$~\cite{Grahamnew,GW16a} (and thus their critical points are hypersurfaces that admit
smooth defining functions subject to $S(\bar g,\sigma) = 1 +
\sigma^{d+1} T$~\cite{GW13,GW15,GW161}).
%
%
% \edz{Cite here or earlier??, discuss Willmore surfaces in intro} In $d$ dimensions, generalized Willmore energy functionals exist 
%such that hypersurfaces subject to  Equation~\nn{nottoB} are critical points.
%Both these facts are explained in~\cite{use and ACF tangentially}. 
%%Note that when $\Sigma$ is not orientable, no single defining function for~$\Sigma$ exists, because defining functions must change sign upon traversing $\Sigma$. 
%%Also, 

\medskip

With these preliminaries established we can now state a {\it singular Obata problem}, which is an obvious singular analogue of the problem answered by the Obata Theorem~\ref{OBATA}.

\begin{question}\label{sphqu}
  Let $(S^d,\bar g)$ be the standard round sphere. Does there exist a smooth function $\sigma\in C^\infty S^d$ such that
  $$
  S(\bar g, \sigma)=1\, ,
  $$
  and such that $g=\si^{-2}\bar{g}$ is not Einstein?
  %% $\sigma$ defines a hypersurface 
  %% $\Sigma\hookrightarrow S^d$ 
  %% that is not totally umbilic?
\end{question}

Recall that a hypersurface embedded in a  conformal  manifold
$(M,[\bar g])$ is totally umbilic if the second fundamental form of
$\bar g$ is pure trace, and this condition is independent of the
choice of metric representative $\bar g$.
A necessary condition for the singular metric $g=\sigma^{-2} \bar g$ to be Einstein with negative scalar curvature, away from the hypersurface~$\Sigma$ defined by $\sigma$, is that $\Sigma$ is totally umbilic~\cite{Goal,LeBrunHeaven}. Thus {\em if we fix the hypersurface } $\Sigma$, the trace-free second fundamental form obstructs the existence of Einstein solutions to~\nn{sYam} with $\mathcal{Z}(\si)=\Sigma$. 
%%Hence this question is an analog of the one answered by the Obata
%%Theorem~\ref{OBATA}.
%
%
%
%Moreover,  a necessary condition for the singular metric $g=\sigma^{-2} \bar g$ to be Einstein with negative scalar curvature away from the hypersurface~$\Sigma$ defined by $\sigma$ is that $\Sigma$ is totally umbilic~\cite{LeB-Heaven,Goal}. Hence this question is an analog of the one answered by the Obata Theorem~\ref{OBATA}. 
The Clifford torus of the next example gives an interesting instance of a Willmore surface embedded in $S^3$ that is not umbilic.

\newcommand{\z}{\mbox{\it\textctz}}

\begin{example}\label{Cliff}
To construct the Clifford torus, first note that $S^3$ equi\-pped~with its standard conformally flat class of metrics may be realized as the space of lightlike lines ${\mathcal N}_+$ in ${\mathbb R}^5\ni(x,y,X,Y,\z)$ (or alternatively, the ray projectivization of the future null cone) defined by 
$$
x^2+y^2+X^2+Y^2=\z^2\, .
$$
The pullback of the Minkowski metric $$\tilde g=dx^2+dy^2+dX^2+dY^2-d\z^2$$
 to the $\z=1$ section of ${\mathcal N}_+$ is the standard round sphere metric $\bar g$ on $S^3$, while conformally related metrics are obtained by  changing this choice of section. 
Calling
$$
r^2=x^2+y^2\mbox{ and } R^2=X^2+Y^2\, ,
$$
the surface $\Sigma=\{r=\frac1{\sqrt{2}}=R,\z=1\}\hookrightarrow S^3$
is Willmore; indeed $\Sigma$ is the Clifford torus. That~$\Sigma$ is Willmore may be easily  verified using Theorem~\ref{obstrThm}.
For that, note
that the  defining function\begin{equation}\label{firstsigma}
\sigma = \Big[\frac{r-R}{\sqrt{2}}\, \Big(
1-\frac{2}{3} \, \Big(\frac{r-R}{\sqrt{2} z}\Big)^2\, 
\Big)\Big]^\star\in C^\infty S^3
\end{equation}
obeys $S(\bar g,\sigma)=1+\sigma^4 T$ with $T$ smooth; in the above $\star$ denotes the pullback from $C^\infty{\mathbb R}^5\to C^\infty S^3$.
To check this computation, note that if $\sigma=\tilde\sigma^\star$ for some $\tilde\sigma\in C^\infty {\mathbb R}^5$ of homogeneity one, then using the ambient formula for  the 
Thomas $D$-operator, the $S$-curvature enjoys an ambient formula
$$
S(\bar g,\sigma)=\Big[
\big|d\tilde\sigma\big|^2_{\tilde g}-\tfrac23\, \tilde\sigma\hh \Delta^{\tilde g} \tilde \sigma
\Big]^\star\, ;
$$
see for example~\cite{CapGoamb}.
A cartoon depiction of this model for the Clifford torus is drawn below.
\begin{center}
\begin{picture}(200,200)
\thicklines
\put(10,44){\rotatebox{2.5}{\includegraphics[width=6cm,height=4cm]{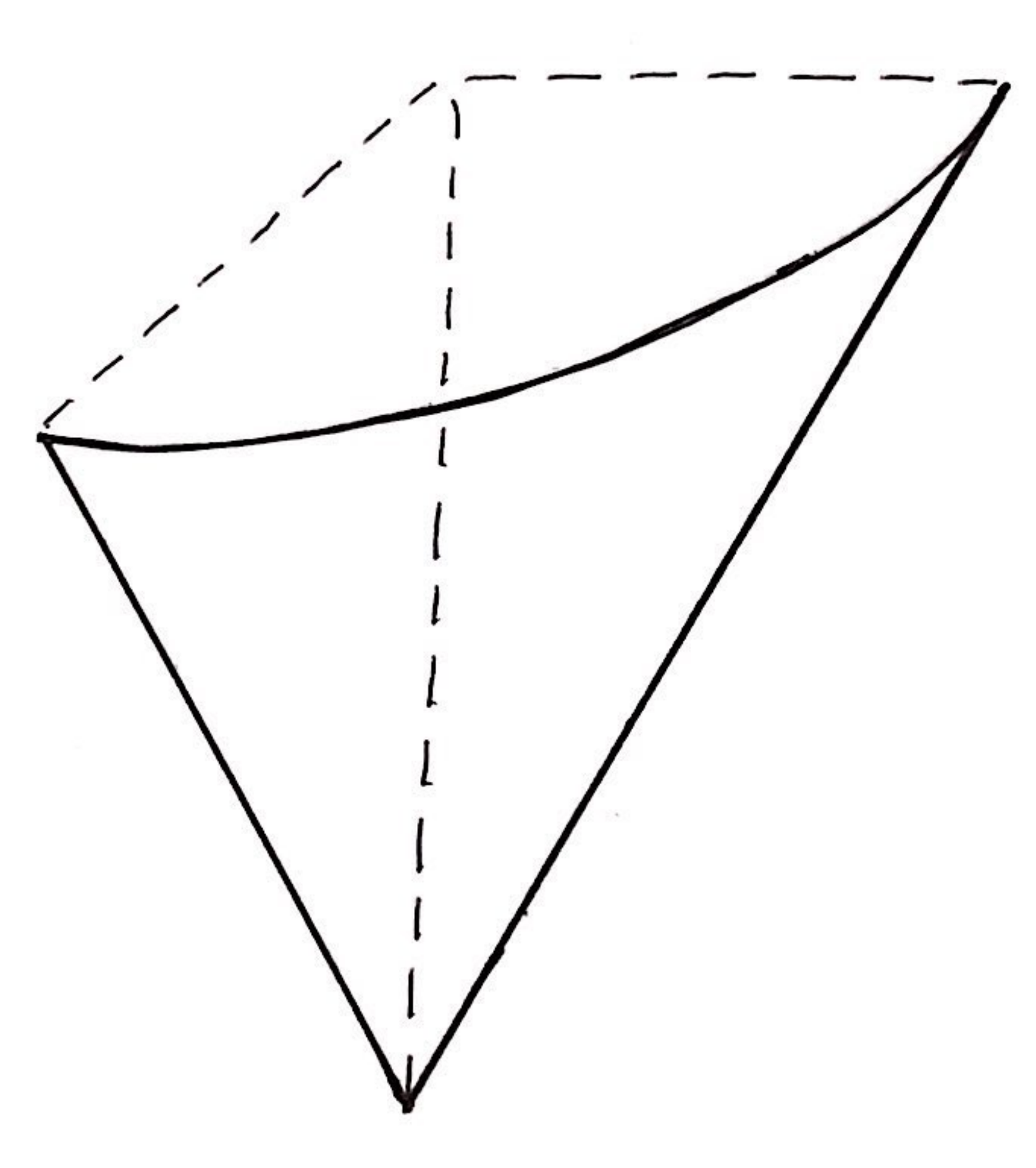}}}
\put(85,160){\vector(0,1){15}}
\put(87,173){$\z$}
\put(64,38){\vector(-1,-1){12}}
\put(44,27){$r$}
\put(100,48){\vector(1,0){15}}
\put(117,45){$R$}
\put(105,124){$\bm \times$}
\put(100,132){$\Sigma$}
\put(133,90){${\mathcal N}_+$}
\put(61,106){$S^3$}
\put(12,165){${\mathbb R}^5$}
\end{picture}
\end{center}
The angle coordinates are suppressed in this picture, so each point (save for when a radius  $r$ or $R$ vanishes) represents a   torus $T^2$.

By using polar coordinates $(r,\theta)$ and $(R,\Theta)$ for the $xy$ and $XY$-planes of the $\z=1$ hyperplane, calling
$$
\tan\tau=\left[\frac{r-R}{r+R}\right]^\star \in [-1,1] \mbox{ where } -\pi/4\leq \tau\leq \pi/4\, ,
$$
the round~$S^3$ metric  in coordinates $(\tau,\theta, \Theta)$  becomes
%$$
%ds^2 = \frac{dt^2}{(1+t^2)^2}+\frac12 \frac{(t+1)^2 d\theta^2 + (t-1)^2 d\Theta^2 }{1+t^2}\, .
%$$
%YIKES, IF WE CALL $t=\tan(\tau/2)$ where $\tau\in [-\pi/2,\pi/2]$
%then
\begin{equation*}\label{nicecoords}
ds^2 =  d\tau^2+\frac12 (d\theta^2 + d\Theta^2) + \frac 12 \sin 2\tau\,  (d \theta^2 - d\Theta^2)\, .
\end{equation*}
Thus, away from circles at $\tau=\pm \pi/4$ (so $r^\star=0$ or $R^\star=0$),
the sphere~$S^3$ is foliated by  constant $\tau$,  torii.  The two circles correspond to a radius of a  torus cycle degenerating to zero.
The 
Clifford torus is the zero mean curvature torus  at $\tau=0$
and has flat torus metric 
$$
ds^2_\Sigma=\frac12(d\theta^2+d\Theta^2)\, .
$$
The second fundamental form for the embedding $\Sigma\hookrightarrow S^3$ in these coordinates is 
$$
\II=\frac12 (d\theta^2 - d\Theta^2)\, .
$$
This is trace-free, {\it i.e.} $\II=\IIo$,  so $\Sigma$ is a minimal surface, but the embedding $\Sigma\hookrightarrow (S^3,[\bar g])$ is not umbilic.
\end{example}

Question~\ref{sphqu} leads to other natural questions, all of which   have
natural extensions to generally curved manifolds. To develop these, we
first present key elements of the theory of conformal hypersurface
embeddings.

\section{Tractor calculus for hypersurface embeddings}\label{tractors}

We say that  a conformal manifold $(M,\cc)$, 
is {\it almost Einstein} (AE) if there exists a metric $\bar g\in \cc$ and smooth function $\sigma$ such that, away from the zero locus ${\mathcal Z}(\sigma)$, the metric 
\begin{equation}\label{ccompact}
g=\frac{\bar g}{\sigma^2}
\end{equation}
has vanishing trace-free Schouten tensor
\begin{equation}\label{AE}
%P(\bar g,\sigma):=
\mathring P^g = 0\, .
\end{equation}
Note that in the special case when $M$ is a $d$-manifold with dimension $(d-1)$ boundary given by ${\mathcal Z}(\sigma)=\partial M$, and in addition $\sigma$  is a defining function for $\partial M$ (so  $\mathrm{d} \sigma \neq
0$ at all points of $\partial M$) then  metrics $g$ obeying~\nn{ccompact} are said to be {\it conformally compact}. When the Einstein condition~\nn{AE} also holds, then such metrics $g$  are called {\it Poincar\'e--Einstein},
or {\it asympotically Poincar\'e--Einstein} when $g$ solves~\nn{AE} asymptotically to the highest order uniquely determined by the  conformal class of metrics on $\partial M$ determined by $\bar g$, see~\cite{FG} for details.

Returning to the more general AE setting, notice that  because $g$ only depends on the ratio $\bar g/\sigma^{2}$, and thus equivalently only on equivalence classes of metric function pairs $[\bar g,\sigma]=[\Omega^2 g,\Omega \sigma]$ with $0<\Omega\in C^\infty M$,  the AE condition is better stated in terms of 
conformal densities defined as follows: A {\it weight $w\in {\mathbb C}$ conformal density} is a  section of the line bundle $$\ce M[w]:=\big((\wedge^d M)^2\big)^{\frac{w}{2d}}\, .$$
The function $\sigma$ defines a section $\bm \sigma$ of $\ce M[1]$ and the volume form $\omega_{\bar g}$  one of $\ce M[d]$. These are related by $\bm \sigma = \sigma\,  (\omega_{\bar g})^{\frac1d}$, which we may also denote by the pair $[\bar g,\sigma]$.
Moreover, tautologically,  $\bar g$ defines  a section  $\bm g\in  \Gamma(\odot^2 T^*M\otimes \ce M[2])$ by $\bm g:= \bar g\,  (\omega_{\bar g})^\frac2d$; we call $\bm g$ the {\it conformal metric}.  
%We may also view the conformal class of metrics $\bm c$ as 
%ray subbundle of $\odot^2 T^*M\otimes\ce M[2]$. 
In addition, where this is defined,  
$g = \bm \sigma^{-2} \bm g$. In general, a {\it true scale} $0<\bm \tau\in \Gamma(\ce M[1])$ canonically defines a metric $g_{\bm \tau}:=\bm \tau^{-2} \bm g$.

The AE Equation~\nn{AE} is overdetermined, so it is propitious to study its prolongation to a triple of sections
$$
(\bm \sigma, \bm n, \bm \rho)\in \Gamma(\ce M[1]\oplus T^*M[1]\oplus \ce M[-1])\, ,
$$
where for any vector bundle $B$ over $M$, we define $B[w]:=B\otimes \ce M[w]$. Then, given any true scale $\bm \tau$, the conformal manifold $(M,\cc)$ is AE iff~\cite{BEG,Goal}
\begin{eqnarray*}
&\nabla^{\bm \tau} \bm \sigma - \bm n = 0\, ,\\[1mm]
&\nabla^{\bm \tau} \bm n + \bm \sigma P^{g_{\bm \tau}} + \bm \rho\,  \bm g = 0\, ,\\[1mm]
&\nabla^{\bm \tau} \bm \rho - P^{g_{\bm \tau}}(\bm n,\pdot)=0\, .
\end{eqnarray*}
In the above, for example,  $\nabla^{\bm \tau} \bm \sigma=[g_{\bm \tau},\nabla^{g_{\bm \tau}}\sigma_{\bm \tau}]\in \Gamma(T^*M[1])$ and the inverse of the conformal metric $\bm g$ is used to contract $\bm n$ with the Schouten tensor; the remaining new notations should then be self explanatory.
The choice of~$\bm \tau$ in the above condition is irrelevant, since a different choice gives an independent linear combination of the three stated equations.
In particular, away from ${\mathcal Z}(\bm \sigma)$, the choice $\bm \tau= \bm \sigma$ gives that  $P^g$ is proportional to the metric and hence that the trace-free Schouten tensor~$\mathring P^g=0$.

In fact, the above system can be re-expressed as a linear connection acting on a section of a suitable, conformally invariant bundle. This bundle is the {\it tractor bundle}
$$
\ct M = \ce M[1]\lpl
T^*M[1]\lpl \ce M[-1]\, ,
$$
and its sections are called (standard) tractors.
The semi-direct sum notation~$\lpl$ indicates that the tractor bundle is a disjoint union of Whitney sum vector bundles $\ce M[1]\oplus T^*M[1]\oplus \ce M[-1]$ indexed by $\bar g\in \cc$ (so equivalently true scales), quotiented by the equivalence  on sections with indices labeled by true scales $\bm \tau$ and $\bm \tau'$, given by
\begin{equation}\label{letschange}
(\bm v^+, \bm v, \bm v^-)_{{\bm \tau}}\sim
\Big(\bm v^+, \bm v+\Upsilon\bm v^+, \bm v^- 
-  {\bm g}^{-1}\!{}
\big( 
\Upsilon,{\bm v}
+\tfrac12 \bm \tau^2  \Upsilon {\bm v}^-\big)\Big)_{{\bm \tau}'}\, ,
\end{equation}
where $\Upsilon :=  \ext\log(\bm \tau'/\bm \tau)$. Observe that the above formula implies that the first non-zero entry of a 
tractor is conformally invariant, this is called the {\it projecting part}. The appropriate linear connection is the {\it tractor connection} $\nabla^\ct :\Gamma(\ct M)\to \Gamma(\ct M\otimes T^*M)$, defined for a metric labeled by the true scale $\bm \tau$, according to 
\begin{equation}\label{letsconnect}
\nabla^\ct (\bm v^+\!, \bm v, \bm v^-)_{\bm \tau} := 
\big(\nabla^{\bm \tau} \bm v^+\! - \bm v, \nabla^{\bm \tau} \bm v + \bm v^+ P^{g_{\bm \tau}} + \bm v^-\,  \bm g ,\nabla^{\bm \tau} \bm v^- \!- P^{g_{\bm \tau}}(\bm v,\pdot\hh)\big)_{\bm \tau}\, .
\end{equation}
The tractor bundle $\ct M[1]$ enjoys a distinguished section $X$ termed the {\it canonical tractor}, defined in any choice of scale $\bm \tau$ by 
$$
X_{\bm \tau}:=(0, 0, 1)\, .
$$
Note that  ${\bm \tau}^{-1} X$ is a standard tractor and, from Equation~\nn{letsconnect}, we see that the tractor connection obeys a non-degeneracy condition
$$
\btau\hh\nabla^\ct ({\bm \tau}^{-1} X)_{\bm \tau}= (0,\bm g,0)\, .
$$

In the above terms, the AE condition~\nn{AE} 
becomes the parallel condition on sections $I\in \Gamma(\ct M)$~\cite{BEG}
$$
\nabla^\ct I=0\, .
$$
A necessary condition for $I$ to be parallel is that, for any choice of true scale~$\bm \tau$,
this tractor is determined in terms of $\bsigma$ according to
\begin{equation}\label{IamI}
I=\big(
\bsigma, \nabla^{\bm \tau}\bm \sigma,
-\tfrac1d\, \btau^{-2}\hh (\Delta^{g_{\bm \tau}}+J^{g_{\bm \tau}})\bm \sigma
\big)_{\bm \tau}\, .
\end{equation}
Tractors determined in terms of a section $\bm \sigma\in \Gamma(\ce M[1])$ this way are termed {\it scale tractors} (the notation $I_{\bm \sigma}$ will sometimes be employed to make the dependence on $\bm \sigma$ clear). 

The tractor bundle enjoys a conformally invariant tractor metric $h\in \Gamma(\odot^2 \ct M)$ that is preserved by the tractor connection $\nabla^{\ct}$ and given by (suppressing the choice of $\bm \tau$ and recycling the notations for sections  $U,V$ of~$\ct M$ used above)
$$
h(U,V)=\bm u^+ \bm v^- + \bm u^- \bm v^+ + \bm g^{-1} (\bm u, \bm v) \in C^\infty M\, .
$$
Observe that given a scale tractor $I$, the corresponding scale $\bm \sigma$ is given in terms of the tractor metric and canonical tractor by
$$
\bsigma= h(I,X)\, .
$$

Now, if $I$ is the  scale tractor determined by $\bm \sigma\in \Gamma(\ce M[1])$, then
the function
$$
S(\bm \sigma)=h(I,I)=:I^2
$$
equals the $S$-curvature defined in Equation~\nn{scurvy} upon identifying $\bm \sigma$ with the metric-function pair $(\bar g,\sigma)$. Thus the singular Yamabe Equation~\nn{singyam} becomes
$$
I^2=1\, .
$$
Writing this condition out using Equation~\nn{IamI} with $\bm \sigma = [\bar g, \sigma]$ gives
$$
I^2 = |\ext \sigma|_{\bar g}^2 -\frac{2\sigma}d\, (\Delta^{\bar g} \sigma + J^{\bar g}\sigma)=1\, .
$$
Hence, along ${\mathcal Z}(\sigma)$, it follows that $\ext \sigma\neq 0$, so by the implicit function theorem $\Sigma={\mathcal Z}(\sigma)$ is a smoothly embedded hypersurface. The same applies when $I^2>0$. Thus we have
\begin{proposition}[Gover~\cite{Goal}]\label{sep-prop}
Let $(M,\cc)$ be a conformal manifold and~$I$ the scale tractor of $\bm\sigma\in \Gamma(\ce M[1])$ such that
$$
I^2>0\, .
$$
%in  a collar neighborhood of the zero locus $\Sigma$ of $\bm \sigma$. 
Then, if $\Sigma\neq \emptyset$,  
$$
M=M_-\sqcup \Sigma \sqcup M_+\, ,
$$
where $\Sigma$ is a smoothly embedded separating hypersurface. Moreover the complements $M\backslash M_{\mp}$ are conformal compactifications of $$M_{\pm}=\{P\in M|\bm \pm \sigma(P)>0\}\, .$$
\end{proposition}
\noindent
This situation is depicted below.
\begin{center}
\begin{picture}(200,165)
\thicklines
\put(10,44){\rotatebox{2.5}{\includegraphics[width=6cm,height=4cm]{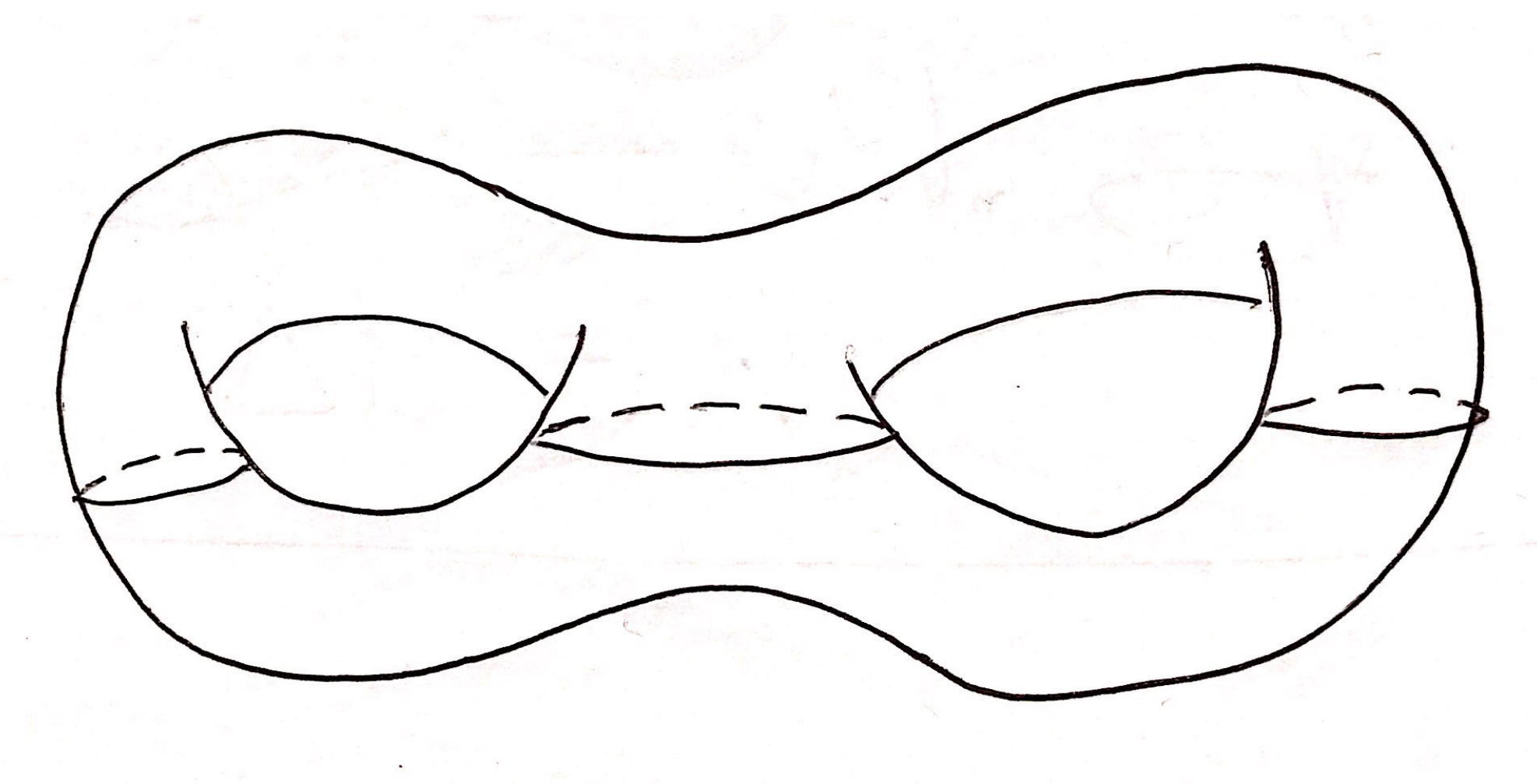}}}
\put(180,100){$\Sigma$}
\put(0,120){$M_+$}
\put(10,60){$M^-$}
\put(80,115){$\sigma>0$}
\put(80,82){$\sigma<0$}
\end{picture}
\end{center}
%\edz{Morse functions, sidedness...?}

The above geometric setup is ideally suited to the study of conformal hypersurface embeddings $\Sigma\hookrightarrow (M,\cc)$. When considering Riemannian hypersurface embeddings $\Sigma\hookrightarrow (M,\bar g)$, a key ingredient is the unit conormal $\hat n$ to $\Sigma$ which may be computed from any defining function $\sigma$ by the formula $\hat n = \ext \sigma'\big|_\Sigma$, where the defining function $\sigma':=\sigma/|\sigma|_{\bar g}$ is a normalized improvement of $\sigma$. 
For conformal embeddings, an analogous, conformally invariant tractor object is the {\it normal tractor} $N\in \Gamma(\ct M|_\Sigma)$, which for a choice of metric $g_{\bm \tau}\in \cc$, is given by
$$
N=(0,\bm {\hat n},-H^{g_{\bm \tau}})_{\bm \tau}\, .
$$
Here we have used that the unit conormal 
rescales as a weight one-density under conformal changes of metric and thus defines $\bm {\hat n}\!\hh\in \!\hh\Gamma(\ct M|_\Sigma)$. Also~$H^{g_{\bm \tau}}$ is the mean curvature of $\Sigma$ computed in the $g_{\bm \tau}$ metric.
 It is not difficult to check that upon replacing the true scale $\bm \tau$ with a new true scale $\bm \tau'$, 
this changes in concordance with Equation~\nn{letschange}. Next we define the conformal analog of a defining function.

\begin{definition}
  A density $\bm \sigma\in \Gamma(\ce M[1])$ satisfying
$$
I^2\big|_\Sigma \neq 0\, , \quad \mbox{where $\Sigma$ is the zero locus of } \si ,
$$
is called a {\it defining density}. 
% \edz{note change here}
%% or equivalently,  whose scale tractor $I$ obeys
%% $$
%% I^2\big|_\Sigma \neq 0\, ,
%% $$
%% is called a {\it defining density}.
\end{definition}\\
\noindent Note that the zero locus of a  defining density subject to $I^2\big|_\Sigma>0$ is a hypersurface.

\medskip 

The next result relates the normal tractor to a defining density for $\Sigma$.
\begin{proposition}[Gover~\cite{Goal,GW161}] \label{IvsN}
Let $\bm \sigma$ be a defining density for $\Sigma$ and 
$$
\bm \sigma' := \bm {\hat \sigma}\, \Big(
1-\tfrac d{4(d-1)} \big(I_{\bm{\hat \sigma}}^2-1\big)
\Big)
\, ,$$
with $\bm {\hat \sigma}:=\bm \sigma/\sqrt{I_{\bm \sigma}^2}$.
Then $$I_{\bm \sigma'}^2 = 1 + \bm \sigma^2 \bm b\, ,$$
for some $\bm b\in \Gamma(\ce M[-2])$, and
$$
I_{\bm \sigma'}\big|_\Sigma = N\, .
$$
\end{proposition}
Observe, that not only does the above proposition give a formula for the normal tractor, it is also the first step in constructing an asymptotic solution to the singular Yamabe problem for $\Sigma$. As we discussed earlier, an important result of Anderson, Chrus\'ciel and Friedrich is that this can be solved to order $d$ in $\bm \sigma$~\cite{Andersson}. In the language of densities, Theorem~\ref{mainonex} is stated as follows.
\begin{theorem} \label{obstrThm}
For any conformally embedded hypersurface  $\Sigma$ $\hookrightarrow$  $(M^d$,  $\cc)$, there exists a defining density $\bm \sigma$, unique up to the addition of smooth terms proportional to~$\bm \sigma^d$,  such that its scale tractor obeys
\begin{equation}\label{normal}
I^2=1+\bm \sigma^d\, \bm B
\end{equation}
for some $\bm B\in \Gamma(\ce M[-d])$.
Moreover $$\bm B_\Sigma=\bm B|_\Sigma$$
is an invariant of the conformal embedding.
\end{theorem} 
The density $\bm B_\Sigma$ is called the {\it obstruction density} because it obstructs smooth solutions to the singular Yamabe problem, where here we mean that the defining density is smooth across $\Sigma$ (when the obstruction density is non-vanishing, one-sided smooth solutions are known to still be possible, but with logarithmic
%\edz{RQ: what's "polylogarithmn''?  Stolen language from, ACF, series including logs....}
behavior around $\Sigma$~\cite{Loewner, Aviles, Mazzeo, Andersson}).  For conformally embedded surfaces, the obstruction density may be obtained as the functional gradient of the Willmore energy functional. Of course, $\bm B_\Sigma$ is the density-valued analog of the higher Willmore invariant introduced earlier.

Another important invariant of conformal hypersurface embeddings is the trace-free part of the second fundamental form $\bm \IIo\in \Gamma(\odot^2 T^*M\otimes\ce M[1]\big|_\Sigma)$. Let us call defining densities that obey Equation~\nn{normal} {\it unit}. Then from Proposition~\ref{IvsN}
and reference~\cite{BEG} (see also \cite{Grant,CG15}) we have the following result. 
\begin{proposition}\label{umb-prop}
Let $\bm \sigma$ be a unit defining density for $\Sigma \hookrightarrow (M,\cc)$. Then,  the projecting part of the restriction of $\nabla^\ct I_{\bm \sigma}$ to $\Sigma$ equals the trace-free second fundamental form $\bm \IIo$.
\end{proposition}
\noindent
Note that  total umbilicity of the boundary of a Poincar\'e--Einstein manifold is a direct corollary of the above proposition.

\medskip
Just as the Gauss formula relates the projection of the ambient Levi-Civita connection to its counterpart intrinsic to the hypersurface $\Sigma$ through the second fundamental form, a similar result holds for the tractor connection. For that, along $\Sigma$  the normal tractor $N$ and tractor metric $h$ give a canonical Whitney direct sum bundle decomposition
$$
\ct M|_\Sigma = \ct M^\perp \oplus \ct M^\parallel\, .
$$
Moreover the rank $d+1$ bundle $\ct M^\perp$ (whose sections $V$ obey $h(V,N)=0$) is isomorphic to the tractor bundle intrinsic to $(\Sigma,\cc_\Sigma)$ where~$\cc_\Sigma$ is the conformal class of metrics on~$\Sigma$ induced from $\cc$ on $M$. From Equation~\nn{letsconnect}, we see that the difference between the projection of the ambient tractor connection to $\ct M^\perp$ and the tractor connection depends on  the difference between their respective Schouten tensors, this is measured by the {\it Fialkow tensor}
$$
{\bm F}\stackrel{\bm \tau}=(P^{g_{\bm \tau}})^\top- P^{g^{\sss\Sigma}_{\bm \tau}}+
H^{g_{\bm \tau}}\,  \bm \IIo +\frac12 
(H^{g_{\bm \tau}})^2
g^{\sss \Sigma}_{\bm \tau}\, \, ,
$$
(an equals sign adorned by $\bm \tau$ is used to indicate a formula computed in a choice of scale $\bm \tau$)
where the ambient dimension $d\geq 4$.
Also, in the above formula, $g_{\bm \tau}^{\sss\Sigma}$ denotes the metric induced by $g_{\bm \tau}$ along~$\Sigma$, the subscript $\top$ projects the tangent bundle of $M$ along~$\Sigma$ to directions orthogonal to the unit conormal, and we have used that this bundle projection gives a rank $(d-1)$ bundle isomorphic to $T\Sigma$.

The conformally invariant, extrinsically coupled tensors $\bm{\hat n}$, $\bm \IIo$, $\bm F$ and~$\bm B_\Sigma$, are key ingredients of a conformal tensor calculus of $\Sigma\hookrightarrow (M,\cc)$. The computation of more interesting quantities
such as higher Willmore energies require
the introduction of a pair of new tractor operators. Firstly, the {\it Thomas $D$-operator} 
gives a mapping 
$$\Gamma(\ct^\Phi M[w])\to  
\Gamma(\ct^\Phi M[w-1]\otimes \ct M)\, ,$$
where $\Phi$ indicates any tensor product of  tractor and conformal density bundles. In a self-explanatory matrix notation (see~\cite{CG15})
$$
D\stackrel{\bm \tau}=
\begin{pmatrix}
w(d+2w-2)\\
(d+2w-2) \nabla^\ct\\
-\btau^{-2}(\Delta^\ct + w J^{g_{\bm \tau}})
\end{pmatrix}\, ,
$$
where $\Delta^\ct$ is the rough Laplace-type operator built in the standard way  from the  metric $g_{\bm \tau}$ and the tractor connection $\nabla^\ct$. Note that the scale tractor $I$ of~$\bm \sigma\in \Gamma(\ce M[1])$ is given by $I=\frac1d D\bm \sigma$.
Contracting the Thomas $D$-operator with the scale tractor $I$ using the tractor metric, gives its {\it Laplace--Robin operator} denoted $I\!\cdot\!D$. Upon restriction to the hypersurface zero locus $\Sigma$ of $\bm \sigma$, this gives a conformally invariant
Robin-type combination of normal and Dirichlet operators due to Cherrier~\cite{cherrier}.
Away from~$\Sigma$, computed in the metric $g_{\bm \sigma}$, the Laplace--Robin operator gives the Laplacian modified by scalar curvature which, when acting on conformal densities at the critical weight $w=1-\frac d2$, recovers the Yamabe operator $-\Delta^{g_{\bm \sigma}} -(1-\frac d2) J^{g_{\bm \sigma}}$. 

Suitable critical powers of the Laplace--Robin operator, upon restriction to~$\Sigma$, give extrinsically coupled Laplacian powers~\cite{GW14,GW161}. These may be viewed as extrinsic generalizations of the conformal Laplacian powers of Graham, Jennes, Mason and Sparling~\cite{GJMS}. 

A log-density is a section of an associated bundle 
induced by the additive representation $y\mapsto y - \log t$ of ${\mathbb R}_+\ni t$ to the bundle of metrics~$\cc$ over $M$, which itself is an ${\mathbb R}_+$ principal bundle. Given a true scale $\bm \tau=[g,\tau]$, the pair $(g,\log \tau)$ determine a log density that we call $\log \bm \tau$. The Laplace--Robin operator can be extended to act on log densities and maps these to standard weight $w=-1$ conformal densities. 
Log densities are essential for writing  compact formul\ae\ for
generalized Willmore energies:
\begin{theorem}
Let $\Sigma\hookrightarrow (M,\cc)$ be a conformally embedded hypersurface and $\bar g = \bm g/\bm \tau^2$ be a metric determined by a true scale $\bm \tau$.  Then
$$
Q^{\bm \tau}:=  (I\!\cdot\! D)^{d-1} \log \bm \tau\big|_\Sigma \in \Gamma(\ce \Sigma [1-d])\, ,
$$ 
where $I$ is the scale tractor of a unit defining density $\bm \sigma$ for $\Sigma$. Moreover, if the singular metric $g=\bm g/\bm \sigma^2$ is asymptotically Poincar\'e--Einstein and~$d$ is odd, then $(-1)^{d-1}((d-2)!!)^{-2} Q_{\bm \tau}$ is the
Branson $Q$-curvature of $(\Sigma,\bm c_{\Sigma})$. Moreover, the integral
$$
{\mathcal A}=\int_\Sigma Q^{\bm \tau}
$$ 
is independent of the choice of $\bar g\in \cc$.
\end{theorem} 
\noindent
In the above we have used that the bundles  $\ce M[w]|_\Sigma$ and $\ce \Sigma[w]$ are isomorphic.
 Results of the above type, where a complicated-to-compute object along a hypersurface is the restriction of a simple bulk quantity, are called {\it holographic formul\ae}. The conformal density $Q^{\bm \tau}$ gives an extrinsic
analog of the  Branson $Q$-curvature 
while its conformally invariant integral $
{\mathcal A}$ is 
a generalized Willmore energy.
In fact, the functional $\frac{(-1)^{d-1}}{(d-1)!\, (d-2)!}\, {\mathcal A}$ is the anomaly coefficient of the log divergence of the renormalized volume
expansion for the singular metric~$g$~\cite{GW16a}.
Moreover, the failure of the  {\it extrinsic $Q$-curvature}  $Q^{\bm \tau}$ to be conformally invariant is controlled by the extrinsic conformal Laplacian powers described above.
Explicit low dimensional  formul\ae\ for the extrinsic $Q$-curvature are given in the following example.

\begin{example}
When $\Sigma\hookrightarrow (M^3,\cc)$ is a conformally embedded surface, 
$$
\int_\Sigma Q^{\bm \tau} = \pi \chi_\Sigma-\frac14  \int_\Sigma {\rm tr}_{\bm g_\Sigma}\hh \bm \IIo^2\, . 
$$
Here $\chi$ denotes the Euler characteristic while the second term is a manifestly conformally invariant bending energy.
For $\Sigma\hookrightarrow (M^4,\cc)$ a closed hypersurface embedded in a conformal four-manifold~\cite{GW161,GGHW}, 
$$
\int_\Sigma Q^{\bm \tau} = \frac 23 \int_\Sigma 
 {\rm tr}_{\bm g_\Sigma} (\hh\bm \IIo\!{}\circ\!{} \bm F)\, .
$$
In the above examples, we have used the conformal metric along $\Sigma$ to construct an endomorphism of $T\Sigma$ from the trace-free second fundamental form and Fialkow tensor; this is the meaning of the symbol ${\rm tr}_{\bm g_\Sigma}$. The extrinsic $Q$-curvature for four manifold embeddings has been computed in~\cite{GrahamRiechert}.
\end{example}

\subsection{The sphere}\label{sphere}
\newcommand{\bsi}{{\bm \sigma}}

Much of the above is nicely illustrated on the sphere equipped with
its standard round conformal structure.  In
this case the conformal structure is  conformally
flat and the  tractor connection~$\nabla^\ct$ has trivial holonomy. Thus
any standard tractor at a point can be extended to a parallel tractor $I$ on the sphere.

In fact this is easily established explicitly, as the conformal $d$-sphere
arises as the ray projectivisation $\mathbb{P}_+$ of the future null cone
$\mathcal{\cN}_+$ of the (Lorentzian signature) Minkowski metric
$\tilde{g}:=\operatorname{diagonal}(-1,1,\cdots ,1)$ on~$\mathbb{R}^{d+2}$ 
(see Example~\ref{Cliff} above). 
Then the tractor connection arises from 
parallel transport in $\mathbb{R}^{d+2}$ viewed as an affine space
(see \cite{CapGoamb,CG15} for more detail). It follows that each
(constant) vector $\tilde{I}$ in $\mathbb{R}^{d+2}$ determines
a corresponding parallel tractor~$I$ on the sphere and vice versa. In
this picture the canonical tractor $X$ 
corresponds to the restriction
to $\cN_+$ of the Euler vector field $\tilde{X}$ of~$\mathbb{R}^{d+2}$. Thus if a parallel tractor~$I$ is timelike,
and normalised to satisfy $I^2=h(I,I)=-1$ say, then $\bsi:=h(I,X)$ has no zero
locus, and the corresponding Einstein metric $g=\bsi^{-2}\bg$ has
scalar curvature~$d(d-1)$. Thus it is the usual round sphere
metric. If on the other hand $I$ is a null
parallel tractor then $\bsi:=h(I,X)$ has an isolated zero
corresponding to where $\tilde{X}$ is parallel to $\tilde{I}$. Since
$I^2=0$ the scalar curvature vanishes and the Einstein metric
$g=\bsi^{-2}\bg$ on the complement is thus Ricci flat. Hence it is
isometric to Euclidean space.  Similarly if a parallel tractor $I$ is
spacelike, satisfying $I^2=1$ then the zero locus $\Sigma$ of
$\bsi:=h(I,X)$ is a hypersurface corresponding to the intersection of
the hyperplane $\tilde{h}(\tilde{I},\tilde{X})=0$ with $\cN_+$, and
the corresponding Einstein metric $g=\bsi^{-2}\bg$ on $S^d\setminus
\Sigma$ has scalar curvature $-d(d-1)$. As it is a conformally flat
Einstein metric it is isometric to the hyperbolic metric. Moreover,
since $I$ is in particular parallel along $\Sigma$, it follows from
Proposition \ref{umb-prop} above that $\Sigma$ is totally umbilic. This, with each of the
connected parts of $S^d\setminus \Sigma$, provides a conformal
compactification  of hyperbolic space which (by stereographic
projection) is conformally equivalent to the usual Poincar\'e-ball.

\begin{center}
	\begin{tikzpicture}[xscale=0.45,yscale=0.45]
		%Cone	
		\draw  (0,3) ellipse (4 and 1);
		\draw (-3.95,2.84) -- (0,-3) -- (3.96,2.86);
		\node at (-2,1.3) {$\mathcal{N}_+$} ;
		\node at (-4,0.6) {$\mathbb{R}^{d+1,1}$};
               % \node at (-4.5,0) {$+$};
%%\node at (-4,0.6) {$\mathbb{R}^{n+2,1}$};

		%Tractor I, subcone and hyperbolic section
		\draw[-latex] (0,-3) --(1.5,-3);
		\node at (1.5,-3.5) {$\tilde{I}$} ;
		 \draw [dotted] (0,-3)  -- (0.99,3.95);
		 \draw [dotted] (0,-3) -- (-0.99,2.03);
		 \draw [dotted] plot[variable=\t,samples=1000,domain=-41.1:49.11] ({tan(\t)+1.33},{8.9*(sec(\t)-1)-0.908});
		\node at (1,4.35) {$\scriptstyle{\tilde{\sigma}=0}$} ;
		\node at (2.5,4.15) {$\scriptstyle{\tilde{\sigma}=1}$} ;
		\node at (3.5,-1.5) {$\tilde{\sigma}= \tilde{I}_A \tilde{X}^A$} ;
		\draw[-latex] (1.5,-0) --(3,0);
		
		%Arrow
		\draw[-latex, thick] (5,0)--(8.5,0);
		\node at (6.7,0.5) {$\mathbb{P}_+$} ;
		
		%Sphere
		\draw  (13,0) circle (3);
		 \draw[dashed,xscale=1,yscale=1,domain=0:3.141,smooth,variable=\t] plot ({0.6*sin(\t r)+13},{3*cos(\t r)});
		 \draw[xscale=1,yscale=1,domain=3.141:6.283,smooth,variable=\t] plot ({0.6*sin(\t r)+13},{3*cos(\t r)});
		\node at (14.75,0) {$\mathbb{H}^{d}$} ;
		\draw[thin] [-stealth] (13,3.7) -- (13,3.06);
		\node at (13.37,4.1) {$S^{d-1}=\partial \mathbb{H}^{d}$};
	\end{tikzpicture}
	\label{Fig:HyperbolicSpaceViaHolonomy}
%\end{center}
\end{center}

\section{Singular Yamabe and  Obata Problems} \label{Qsect}

Let $(M,\cc)$ be a conformal manifold equipped with a smooth and sign
changing solution $\bsi$ of the singular Yamabe problem $S(\bm c, \bm \si)=1$. It
follows that $\Sigma=\mathcal{Z}(\bm \si)\neq \emptyset$, and then from
Proposition \ref{sep-prop} and Theorem~\ref{obstrThm} that~$\Sigma$ is
a smoothly embedded separating hypersurface that satisfies the higher
Willmore equation $\bm B_{\Sigma}=0$.  Thus 
locally, such solutions potentially provide an interesting route to
accessing and studying higher Willmore hypersurfaces \cite{GW13,GW15,GW161}.

The global problem is potentially even more intriguing. The
existence of sign changing solutions $\bm \bsi$ to the singular Yamabe
problem is itself interesting, and then when such exist, the zero locus
$\Sigma$ is a closed embedded higher-Willmore hypersurface. There are
examples with $\Sigma$ totally umbilic and the metric $g=\bsi^{-2}\bg$
Einstein on the complement of~$\Sigma$. We have seen this above on the
sphere in Section \ref{sphere} (following~\cite{Goal}), and there are also examples on
suitable products of the sphere with Einstein manifolds~\cite{GL1,GL2}. This leads to our main questions.

%% One reason that the singular Yamabe problem $S=1$ is interesting is that for
%% any smooth solution $\bsi\in \Gamma (\ce[1])$ where $\bsi$ changes sign it follows at once
%% that 

%% Let $(M,\cc)$ be a closed conformal manifold equipped with a smooth
%% ASC scale \edz{I guess that we need less regularity} $\si\in \Gamma
%% (\ce[1])$ such that $I^2=1$, equivalently the $S$-curvature \edz{notn}
%% satisfies $S=1$. Then (\ref{obstrThm})

%% $d$-odd
\begin{question}\label{one}
  Do there exist closed conformal manifolds $(M,\cc)$ that admit a smooth, sign changing singular Yamabe scale $\bm \si\in \Gamma (\ce[1])$ such the $S$-curvature $S=I^2$
 obeys  $S=1$ but with  $g=\bsi^{-2}\bg$ not Einstein on $M\setminus \mathcal{Z}(\bm \si)$?
\end{question}

In fact the situation would be most interesting if the zero locus is
not totally umbilic. Thus the following is an interesting problem.
\begin{question}\label{two}
  Do there exist closed conformal manifolds $(M,\cc)$ that admit a sign changing singular Yamabe scale  $\bm \si\in \Gamma (\ce[1])$ 
 such the $S$-curvature $S=I^2$
 obeys  $S=1$   
 and  such that $\Sigma= \mathcal{Z}(\bm \si)$ is not totally umbilic?
\end{question}

\noindent
Note that a positive answer to this implies a positive answer to
Question~\ref{one}, because  if $g=\bsi^{-2}\bg$ is Einstein then the scale
tractor $I$ is parallel everywhere, but along $\Sigma$ agrees with the
normal tractor which is thus parallel and so, as discussed earlier,  $\Sigma=
\mathcal{Z}(\bm \si)$ is totally umbilic.

\medskip 
There are refinements of these questions  where we
assume the initial conformal manifold includes an Einstein metric. The first here is a variant of what we
called earlier, the singular Obata problem.
\begin{question}\label{one-E}
  Do there exist closed Einstein manifolds $(M,\bar{g})$ that admit a sign changing singular Yamabe scale $\bm \si\in \Gamma (\ce[1])$ 
   such the $S$-curvature $S=I^2$
 obeys  $S=1$ 
     but with  $g=\bsi^{-2}\bg$ not Einstein on $M\setminus \mathcal{Z}(\bm \si)$?
\end{question}

\noindent
Once again one can ask for a stronger result, as follows.
\begin{question}\label{two-E}
  Do there exist closed Einstein manifolds $(M,\bar{g})$ that admit a
  sign changing singular Yamabe scale $\bm\si\in \Gamma (\ce[1])$ 
    such that the $S$-curvature $S=I^2$
 obeys  $S=1$
   and such
  that $\Sigma= \mathcal{Z}(\bm\si)$ is not totally umbilic?
\end{question}

\noindent
The former of course generalises the Question \ref{sphqu}
for the sphere. 
Questions \ref{two} and \ref{two-E} are particularly interesting because the obstruction density ${\bm B}_\Sigma$ is a higher Willmore invariant~\cite{GW13,GW15,GW161}.

\addtocontents{toc}{\SkipTocEntry}
\section*{Acknowledgements}

A.W.~was also supported by a Simons Foundation Collaboration Grant for Mathematicians ID 317562.
A.W. and A.R.G.
 gratefully acknowledge support from the Royal Society of New Zealand via Marsden Grant 16-UOA-051.

%number,names
\newcommand{\msn}[2]{\href{http://www.ams.org/mathscinet-getitem?mr=#1}{#2}}
%number
\newcommand{\hepth}[1]{\href{http://arxiv.org/abs/hep-th/#1}{arXiv:hep-th/#1}}
%number
\newcommand{\maths}[1]{\href{http://arxiv.org/abs/math/#1}{arXiv:math/#1}}
%number
\newcommand{\mathph}[1]{\href{http://lanl.arxiv.org/abs/math-ph/#1}{arXiv:math-ph/#1}}
\newcommand{\arxiv}[1]{\href{http://lanl.arxiv.org/abs/#1}{arXiv:#1}}


\begin{thebibliography}{KUVV12}

\bibitem{Andersson}
L. Andersson, P. T. Chru\'sciel and H, Friedrich, {\em On the regularity of solutions to the Yamabe equation and the existence of smooth hyperboloidal initial data for Einstein's field equations}, Comm. Math. Phys. {\bf 149} (1992), 587--612.

\bibitem{Aubin}
T. Aubin,  (1976), {\em \'Equations diff\'erentielles non lin\'eaires et probl\`eme de Yamabe concernant la courbure scalaire}, J. Math. Pures Appl., (9) {\bf 55} (1976), 269--296.

\bibitem{Aviles}
P. Aviles and R. C. McOwen, {\em Complete conformal metrics with negative scalar curvature in compact Riemannian manifolds}, Duke. Math. J. {\bf 56} (1988), 395--398.

\bibitem{BEG}
\msn{1322223}{T.~N. Bailey, M.~G. Eastwood and A.~R. Gover,} \textsl{ Thomas's structure
  bundle for conformal, projective and related structures},
\newblock Rocky Mountain J. Math. \textbf{ 24}(4) (1994), 1191--1217.


\bibitem{BQ}
\msn{1316845}{T. P. Branson}, \textsl{Sharp inequalities, the functional determinant, and the complementary series}, Trans. Amer. Math. Soc. {\bf 347}, no. 10, 3671--3742.  (1995).
%
\bibitem{BO}
\msn{1050018}{T. P. Branson and B. \Odane rsted}, 
\textsl{Explicit functional determinants in four dimensions}, Proc. Amer. Math. Soc. {\bf 113} (1991), 669--682.


\bibitem{CapGoamb} \msn{1996768}{A.\ \v Cap, and A.R.\ Gover}, {\em Standard
     tractors and the conformal ambient metric construction},  Ann.\
     Global Anal.\ Geom.\  {\bf 24} (2003), 231--295.


\bibitem{WhatQ}
\msn{2407525}{S.-Y. Chang, M. Eastwood, 
B. \Odane rsted, Paul C. Yang,}\textsl{What is~$Q$-curvature?},
Acta Appl. Math.
{\bf 102}, Issue 2, 119--125 (2008).


\bibitem{cherrier}
\msn{0749522}{P.~Cherrier,} \textsl{ Probl\`emes de {N}eumann non lin\'eaires sur les
  vari\'et\'es riemanniennes},
\newblock J. Funct. Anal. \textbf{ 57}(2), 154--206 (1984).


\bibitem{CG15}
S. Curry and A. R. Gover, \textsl{An introduction to conformal geometry and tractor calculus, with a view to applications in general relativity}, in Asymptotic Analysis in General Relativity, 86--170, Cambridge University Press, 2018, \arxiv{1412.7559}.

 \bibitem{FG} \href{http://www.ams.org/mathscinet-getitem?mr=837196}{C.\ Fefferman and C.R.\ Graham,} {\em Conformal
    invariants} in: The mathematical heritage of \'{E}lie Cartan (Lyon,
  1984).  Ast\' erisque 1985, Numero Hors Serie, 95--116. 

\bibitem{GGHW}
\msn{3903687}{M.~Glaros, A.~R.~Gover, M.~Halbasch and A.~Waldron}, 
\textsl{Singular Yamabe Problem Willmore Energies},
J. Geom. Phys. {\bf 138} 168--193 (2019),  
\newblock \arxiv{1508.01838}.

\bibitem{Goal}
\msn{2587388}{A.~R. Gover,} \textsl{ Almost {E}instein and {P}oincar\'e-{E}instein manifolds
  in {R}iemannian signature},
\newblock J. Geom. Phys. \textbf{ 60}(2) (2010), 182--204, \arxiv{0803.3510}.

\bibitem{GAW}
A. R. Gover, C. Arias and A. Waldron, {\em
Conformal Geometry of Embedded Manifolds with Boundary from Universal Holographic Formul\ae},
 \arxiv{1906.01731}.

\bibitem{GL1}
\msn{2678944}{
A. R. Gover and F. Leitner},
{\em A class of compact Poincaré-Einstein manifolds: properties and construction},
Commun. Contemp. Math. {\bf 12} (2010),  629--659. 

\bibitem{GL2}
\msn{2574315}{A. R. Gover and F. Leitner},
{\em 
A sub-product construction of Poincaré-Einstein metrics},
Internat. J. Math. {\bf 20} (2009), 1263--1287. 


\bibitem{GW14}
\msn{3218267}{A.~R. Gover and A.~Waldron,} \textsl{ Boundary calculus for conformally compact
  manifolds},
\newblock Indiana Univ. Math. J. \textbf{ 63}(1), 119--163 (2014),
\arxiv{1104.2991}.
%


\bibitem{GW13}
A.~R. Gover and A.~Waldron, \textsl{ {Submanifold conformal invariants and a
  boundary Yamabe problem
   {Conference on Geometrical Analysis-Extended Abstract}, CRM Barcelona (2013),
 {\rm arXived as} Generalising the Willmore equation:
  submanifold conformal invariants from a boundary Yamabe problem}},
  \arxiv{1407.6742}.


\bibitem{GW15}
A.~R. Gover and A.~Waldron, \textsl{ {Conformal hypersurface geometry via a
  boundary Loewner-Nirenberg-Yamabe problem}},
  Commun. Anal. Geom. to appear,
\newblock \arxiv{1506.02723}.
%
%
\bibitem{GW16a}
A.~R. Gover and A.~Waldron, \textsl{ {Renormalized Volume}},
Commun. Math. Phys. {\bf 354} (2017) 1205--1244,       
\newblock 
 \arxiv{1603.07367}.
%
\bibitem{GW161}
A.~R. Gover and A.~Waldron, \textsl{ 
A calculus for conformal hypersurfaces and new higher Willmore energy functionals}, Adv. Geom. in press, \arxiv{1611.04055}.









\bibitem{Grahamnew} \msn{3601568}{C. R. Graham},
\textsl{Volume renormalization for singular Yamabe metrics},
Proc. Amer. Math. Soc. {\bf 145}, 1781--1792 (2017), \arxiv{1606.00069}.


\bibitem{GJMS}
\msn{1190438}{C.R.~Graham, R.~Jenne, Ralph, L.~Mason and G.~Sparling},
{\em Conformally invariant powers of the Laplacian. I. Existence}, 
J. London Math. Soc. (2) {\bf 46} (1992), 557--565.
%

\bibitem{GrahamRiechert}
C.R. Graham and N. Riechert, \textsl{
Higher-dimensional Willmore energies via minimal submanifold asymptotics},
\arxiv{1704.03852}.

\bibitem{Grant}
D.~Grant,
\newblock
  \href{http://www.math.auckland.ac.nz/mathwiki/images/5/51/GrantMSc.pdf}{\it A
  conformally invariant third order Neumann-type operator for hypersurfaces},
\newblock Master's thesis, University of Auckland, New Zealand, 2003.



\bibitem{LeBrunHeaven}
\msn{0652038}{C. R. LeBrun,} {\em~$\mathscr{H}$-Space with a Cosmological Constant},
Proc. R. Soc. Lond.  A {\bf  380}, (1982) 171--185. 


\bibitem{Loewner} C. Loewner and L. Nirenberg, {\it Partial differential equations invariant under conformal or projective transformations}, in Contributions to analysis (a collection of papers dedicated to Lipman Bers), 245--272, Academic Press, New York, 1974.

\bibitem{Marques}
\msn{3152944}{F.C. Marques and A. Neves}, {\em Min-max theory and the {W}illmore conjecture},
 Ann. of Math. (2),
   {\bf 179},
    683--782 (2014).

\bibitem{Mazzeo}
R. Mazzeo, {\it Regularity for the singular Yamabe problem}, Indiana Univ. Math. J. {\bf 40} (1991), 1277--1299.




\bibitem{Obata}
M. Obata, {\em The conjectures on conformal transformations of Riemannian manifolds}, J. Diff. Geom. {\bf 6} (1972), 247--258.



\bibitem{Schoen}
R. Schoen, {\em Variational Theory for the Total Scalar Curvature Functional
for Riemannian Metrics and Related Topics}, in Topics in Variational Calculus, Lecture Notes in Mathematics {\bf 1365}, 120--154, Springer 1980.

\bibitem{Schoen1} R. Schoen
{\em Conformal deformation of a Riemannian metric to constant scalar curvature}, J. Diff. Geom., {\bf 20} (1984), 479--495.

\bibitem{Trudinger}
N.S. Trudinger,  {\em Remarks concerning the conformal deformation of Riemannian structures on compact manifolds}, Ann. Scuola Norm. Sup. Pisa (3), {\bf 22} (1968), 265--274.

\bibitem{Willmore}
\msn{0202066}{T.J. Willmore},
    {\em Note on embedded surfaces},
  An. \c Sti. Univ. ``Al. I. Cuza'' Ia\c si Sec\c t. I a Mat. (N.S.) {\bf 11B}, 493--496 (1965).
  
\bibitem{WillmoreB} T.J.  Willmore, T. J. (1992), {\em A survey on Willmore immersions} in Geometry and Topology of Submanifolds IV, 11-16,  World Scientific 1991.

\bibitem{Yamabe}
H. Yamabe (1960), {\em On a deformation of Riemannian structures on compact manifolds}, Osaka J.  Math., 12 (1960), 21--37.

%
%%[{\SSmidge}A{\SSmidge}{\SSmidge}G{\SSmidge}M{\SSmidge}O$^{\!2}${\SSmidge}0{\SSmidge}0{\SSmidge}]
%
%%\bibitem[{\SSmidge}A{\SSmidge}{\SSmidge}G{\SSmidge}M{\SSmidge}O{\SSmidge}00]{AdSCFTreview}
%%\msn{1743597}{O.\ Aharony, S. S.\ Gubser, J.M.\ Maldacena, H.~Ooguri and Y.~Oz,}
%%  \textsl{Large~$N$ field theories, string theory and gravity},
%%  Phys.\ Rept.\  {\bf 323}, 183--386 (2000),
%%  \href{http://arxiv.org/abs/hep-th/9905111}{arXiv:hep-th/9905111}.
%
%\bibitem[And01]{Anderson}
%\msn{1825268}{M. Anderson}, \textsl{$L^2$ curvature and volume renormalization of the AHE metrics on 4-manifolds}, Math. Res. Lett. {\bf 8} (2001) 171--188,
%\maths{0011051}.
%
%
%\bibitem[ACF92]{ACF}
%\msn{1186044}{L.~Andersson, P.~T. Chru{\'s}ciel and H.~Friedrich,} \textsl{ On the regularity
%  of solutions to the {Y}amabe equation and the existence of smooth
%  hyperboloidal initial data for {E}instein's field equations},
%\newblock Comm. Math. Phys. \textbf{ 149}(3), 587--612 (1992),  \arxiv{0802.2250}.
%
%
%%\bibitem[AGS14]{Astaneh}
%%A.~F. Astaneh, G.~Gibbons and S.~N. Solodukhin, \textsl{ {What surface
%%  maximizes entanglement entropy?}},
%%\newblock Phys. Rev. \textbf{ D90}(8), 085021--085031 (2014), \arxiv{1407.4719}.
%
%
%\bibitem[AM10]{Alexakis}
%\msn{2653898}{S.~Alexakis and R.~Mazzeo,} \textsl{Renormalized area and properly embedded
%  minimal surfaces in hyperbolic 3-manifolds},
%\newblock Comm. Math. Phys. \textbf{ 297}(3), 621--651 (2010), \arxiv{0802.2250}.
%
%
%\bibitem[Au76]{Aubin}
%\msn{0431287}{T. Aubin}, \textsl{\'Equations diff\'erentielles non lin\'eaires et probl\`eme de Yamabe concernant la courbure scalaire}, J. Math. Pures Appl. {\bf 55}, 269--296 (1976).
%
%

%
%
%\bibitem[BJ10]{BaumJuhl}
%H. Baum and A. Juhl, \textsl{Conformal Differential Geometry: Q-Curvature and Conformal Holonomy},
%Birkh\"auser, 2010.
%
%%\bibitem[BC16]{complexity2}
%%O. Ben-Ami and D. Carmi,
%%\textsl{On Volumes of Subregions in Holography and Complexity},
%%\arxiv{1609.02514}.
%
%
%
%
%\bibitem[B95]{BQ}
%\msn{1316845}{T. P. Branson}, \textsl{Sharp inequalities, the functional determinant, and the complementary series}, Trans. Amer. Math. Soc. {\bf 347}, no. 10, 3671--3742.  (1995).
%
%%\bibitem[BO91]{BO}
%%\msn{1050018}{T. P. Branson and B. \Odane rsted}, 
%%\textsl{Explicit functional determinants in four dimensions}, Proc. Amer. Math. Soc. {\bf 113} (1991), 669--682.
%%
%%\bibitem[BG01]{BrGoCNV}
%%\msn{1867890}{T. P.~Branson and A.~R. Gover,} \textsl{ Conformally invariant non-local
%%  operators},
%%\newblock Pacific J. Math. \textbf{ 201}(1), 19--60 (2001).
%
%
%%\bibitem[B\!\hh RSSZ1\!\hh 5]{complexity1}
%%A. R. Brown, D. A. Roberts, L. Susskind, B. Swingle, Y. Zhao
%%\textsl{Holographic Complexity Equals Bulk Action?}
%%Phys. Rev. Lett. {\bf 116} (2016), 191301,
%%\arxiv{1509.07876}.
%
%%\bibitem[CSS01]{CSS} 
%%\msn{1847589}{A. \v{C}ap, J. Slov{\'a}k  and V. Sou\v{c}ek,} \textsl{Bernstein-Gelfand-Gelfand sequences}, Ann. of Math. {\bf 154}, 97--113 (2001), \maths{0001164}.
%%
%%
%%\bibitem[CEOY08]{WhatQ}
%%\msn{2407525}{S.-Y. Chang, M. Eastwood, 
%%B. \Odane rsted, Paul C. Yang,}
%%\textsl{What is~$Q$-curvature?},
%%Acta Appl. Math.
%%{\bf 102}, Issue 2, 119--125 (2008).
%
%% \bibitem[CG00]{CapGoirred}\msn{2247867}{A.\ \v Cap, and A.R.\ Gover,} {\em Tractor
%%       bundles for irreducible parabolic geometries.  Global analysis
%%       and harmonic analysis}, S\'emin. Congr. {\bf 4}, 129, Soc. Math.
%%     France 2000. 
%%
%
%
%
%
%%\bibitem[CG02]{CapGoTAMS} \msn{1873017}{A.\ \v Cap, and A.R.\ Gover,} {\em Tractor calculi
%%    for parabolic geometries}, Trans.\ Amer.\ Math.\ Soc.\ {\bf 354}
%%  (2002), 1511--1548. 
%%
%%
%\bibitem[CG03]{CapGoamb} \msn{1996768}{A.\ \v Cap, and A.R.\ Gover}, {\em Standard
%     tractors and the conformal ambient metric construction},  Ann.\
%     Global Anal.\ Geom.\  {\bf 24} (2003), 231--295.
%
%\bibitem[Cas18]{Case}
%\msn{3776023}{J.S. Case}, 
%\textsl{Boundary operators associated with the Paneitz operator}, Indiana Univ. Math. J. {\bf 67} 293--327 (2018),
%\arxiv{1509.08342 }.
%
%\bibitem[Cha97]{Chang-Qing}
%\msn{1454485}{S.-Y. Alice Chang and J. Qing},
%\textsl{The zeta functional determinants on manifolds with boundary.~I. The formula}, J. Funct. Anal. {\bf 147} (1997),  327--362.
%
%\bibitem[Cha08]{whatQ}
%\msn{2407525}{S.-Y. Alice Chang, Michael Eastwood, Bent \O rsted, Paul C. Yang},
%\textsl{What is Q-Curvature?}
%Acta Applicandae Mathematicae
%{\bf 102}, 119--126 (2008).
%
%\bibitem[Cha081]{Chang}
%\msn{2336463}{S.-Y. Alice Chang, J.J. Qing and P.Yang},
%\textsl{On the renormalized volumes for conformally compact Einstein manifolds}, J. Math. Sci. (N. Y.) {\bf 149} (2008), 1755--1769, \maths{0512376}.
%
%\bibitem[Che84]{cherrier}
%\msn{0749522}{P.~Cherrier,} \textsl{ Probl\`emes de {N}eumann non lin\'eaires sur les
%  vari\'et\'es riemanniennes},
%\newblock J. Funct. Anal. \textbf{ 57}(2), 154--206 (1984).
%
%%\bibitem[CG15]{CG15}
%%S. Curry and A. R. Gover, \textsl{An introduction to conformal geometry and tractor calculus, with a view to applications in general relativity}, \arxiv{1412.7559}.
%
%%\bibitem[ES97]{EastwoodSlovak}
%%\msn{1483772}{M.~Eastwood and J.~Slov{\'a}k,} \textsl{ Semiholonomic {V}erma modules},
%%\newblock J. Algebra \textbf{ 197}(2), 424--448 (1997).
%%
%%
%%\bibitem[EW14]{Engelhardt} 
%%  N.~Engelhardt and A.~C.~Wall,
%%  \textsl{Quantum Extremal Surfaces: Holographic Entanglement Entropy beyond the Classical Regime},
%%  JHEP {\bf 1501}, 073--098 (2015), \arxiv{1408.3203}.
% 
% \bibitem[DM08]{MalDj}    \msn{2456884}{Z. Djadli and A. Malchiodi}, \textsl{Existence of conformal metrics with constant Q-curvature}, Ann. of Math. {\bf 168}, 813--858 (2008), \maths{0410141}.
% 
% 
% 
%\bibitem[Fi44]{Fialkow}
%\msn{0011023}{A.  
%Fialkow},\textsl{
%Conformal differential geometry of a subspace},
%Trans. Amer. Math. Soc. {\bf 56}, 309--433 (1944). 
%
% \bibitem[FG84]{FGast} \href{http://www.ams.org/mathscinet-getitem?mr=837196}{C.\ Fefferman, and C.R.\ Graham,} {\em Conformal
%    invariants} in: The mathematical heritage of \'{E}lie Cartan (Lyon,
%  1984).  Ast\' erisque 1985, Numero Hors Serie, 95--116. 
%  
%\bibitem[FG07]{FGrnew} \msn{2858236}{C.\ Fefferman, and C.R.\ Graham}, {\em The Ambient Metric},
%  Annals of Mathematics Studies, {\bf 178}, Princeton University Press, \arxiv{0710.0919}. 
%  %Cited on pages
%
%
%\bibitem[FG02]{FGQ} \msn{1909634}{C. Fefferman and C. R. Graham,} \textsl{$Q$-curvature and Poincar\'e metrics}, Math. Res. Lett. {\bf 9},  no. 2-3, 139--151 (2002), \maths{0110271}.
%
%\bibitem[FH03]{Fefferman-Hirachi}
%\msn{2025058}{C. Fefferman, and K. Hirachi}, \textsl{Ambient metric construction of Q-curvature in conformal and CR
%geometries}, Math. Res. Lett. {\bf 10}, 819--832 (2003), \maths{0303184}.
%
%\bibitem[FG12]{FGQJuhl}
%\msn{3073887}{C. Fefferman and C.R.~Graham}, {\it Juhl's Formulae for GJMS Operators and Q-Curvatures},
%J. Amer. Math. Soc. {\bf 26},
%1191--1207 (2013), 
%\arxiv{1203.0360}.
%
%
%%\bibitem[GJMS92]{GJMS}
%%\msn{1190438}{C.R.~Graham, R.~Jenne, Ralph, L.~Mason and G.~Sparling},
%%{\em Conformally invariant powers of the Laplacian. I. Existence}, 
%%J. London Math. Soc. (2) {\bf 46} (1992), 557--565.
%
%
%\bibitem[{\SSmidge}G{\SSmidge}G{\SSmidge}HW{\SSmidge}1{\SSmidge}5{\SSmidge}]{GGHW15}
%{\SSmidge}
%\msn{3903687}{M.~Glaros, A.~R.~Gover, M.~Halbasch and A.~Waldron}, 
%\textsl{Singular Yamabe Problem Willmore Energies},
%J. Geom. Phys. {\bf 138} 168--193 (2019),  
%\newblock \arxiv{1508.01838}.
%
%
%
%
%\bibitem[Gov07]{GoSigma}
%\msn{2366922}{A.~R. Gover,} \textsl{ Conformal {D}irichlet-{N}eumann maps and
%  {P}oincar\'e-{E}instein manifolds},
%\newblock SIGMA Symmetry Integrability Geom. Methods Appl. \textbf{ 3}, 
%  100--121 (2007), \arxiv{0710.2585}.
%
%\bibitem[Gov10]{Goal}
%\msn{2587388}{A.~R. Gover,} \textsl{ Almost {E}instein and {P}oincar\'e-{E}instein manifolds
%  in {R}iemannian signature},
%\newblock J. Geom. Phys. \textbf{ 60}(2), 182--204 (2010), \arxiv{0803.3510}.
%
%\bibitem[GP03]{GoPetCMP} \msn{1822358}{A.R.\ Gover, and L.\ Peterson,} {\em Conformally
%    invariant powers of the Laplacian,~$Q$-curvature, and tractor
%    calculus}, Comm.\ Math.\ Phys.\  {\bf 235},  339--378  (2003), \mathph{0201030}. 
%    Cited on pages
%
%\bibitem[GP18]{GPnew}
%A. R. Gover and  L. Peterson,
%\textsl{Conformal boundary operators, T-curvatures, and conformal fractional Laplacians of odd order},
%\arxiv{1802.08366}.
%
%\bibitem[GLW15]{GLW}
%\msn{3338300}{A.~R. Gover, E.~Latini and A.~Waldron,} \textsl{ Poincar\'e-{E}instein
%  holography for forms via conformal geometry in the bulk},
%\newblock Mem. Amer. Math. Soc. \textbf{ 235}(1106), vi+95 (2015),
%\arxiv{1205.3489}.
%
%\bibitem[GSW08]{GoverS} 
%  \msn{2502272}{A.~R.~Gover, A.~Shaukat and A.~Waldron,}
%  \textsl{Tractors, Mass and Weyl Invariance},
%  Nucl.\ Phys.\ B {\bf 812}, 424--455 (2009), \arxiv{0812.3364};
%  \textsl{Weyl Invariance and the Origins of Mass,}
%  Phys.\ Lett.\ B {\bf 675}, 93--97 (2009), \arxiv{0810.2867}.
%%
%%
%%\bibitem[GSS08]{GSS08}
%%\msn{2372762}{A. R. Gover, P. Somberg and V. Sou\v{c}ek,} \textsl{
%%Yang-Mills detour complexes and conformal geometry},
%%Comm. Math. Phys. {\bf 278} 307--327 (2008), \maths{0606401}.
%
%
%\bibitem[GW13]{CRMouncementCRM}
%A.~R. Gover and A.~Waldron, \textsl{ {Submanifold conformal invariants and a
%  boundary Yamabe problem
%   {Conference on Geometrical Analysis-Extended Abstract}, CRM Barcelona (2013),
% {\rm arXived as} Generalising the Willmore equation:
%  submanifold conformal invariants from a boundary Yamabe problem}},
%  \arxiv{1407.6742}.
%
%\bibitem[GW14]{GW}
%\msn{3218267}{A.~R. Gover and A.~Waldron,} \textsl{ Boundary calculus for conformally compact
%  manifolds},
%\newblock Indiana Univ. Math. J. \textbf{ 63}(1), 119--163 (2014),
%\arxiv{1104.2991}.
%
%\bibitem[GW15]{GW15}
%A.~R. Gover and A.~Waldron, \textsl{ {Conformal hypersurface geometry via a
%  boundary Loewner-Nirenberg-Yamabe problem}},
%  Commun. Anal. Geom. to appear,
%\newblock \arxiv{1506.02723}.
%
%
%\bibitem[GW16a]{GWvol}
%A.~R. Gover and A.~Waldron, \textsl{ {Renormalized Volume}},
%Commun. Math. Phys. {\bf 354} (2017) 1205--1244,       
%\newblock 
% \arxiv{1603.07367}.
%
%\bibitem[GW16b]{GW161}
%A.~R. Gover and A.~Waldron, \textsl{ 
%A calculus for conformal hypersurfaces and new higher Willmore energy functionals}, Adv. Geom. in press, \arxiv{1611.04055}.
%
%\bibitem[GW19]{GWvolII}
%A.~R. Gover and A.~Waldron, \textsl{ {Renormalized Volumes with Boundary}}, 
%Commun. Contemp. Math. {\bf 21},
%1850030--1850061 (2019), \arxiv{1611.08345 }.
%
%
%\bibitem[Gra00]{Gra00} \msn{1758076}{C. R. Graham,}
%\textsl{Volume and area renormalizations for conformally compact Einstein metrics},
%Proceedings of the 19th Winter School ``Geometry and Physics'' (Srn\'i, 1999). 
%Rend. Circ. Mat. Palermo (2) Suppl. No. {\bf 63}, 31--42 (2000), \maths{9909042}.
%
%
%
%
%
%
%\bibitem[Gra16]{Grahamnew} \msn{3601568}{C. R. Graham},
%\textsl{Volume renormalization for singular Yamabe metrics},
%Proc. Amer. Math. Soc. {\bf 145}, 1781--1792 (2017), \arxiv{1606.00069}.
%
%%\bibitem[GH05]{GraHi}
%%\msn{2160867}{C.~R.~Graham and K.~Hirachi,}
%%\textsl{The ambient obstruction tensor and~$Q$-curvature}, in  \textsl{AdS/CFT correspondence: Einstein metrics and their conformal boundaries}, 
%%IRMA Lect. Math. Theor. Phys. {\bf 8}, 59--71, European Math. Society, Z\"urich, 2005, \maths{0405068}.
%
%\bibitem[GJMS92]{GJMS}
%\msn{1190438}{C.~R. Graham, R.~Jenne, L.~J. Mason and G.~A.~J. Sparling,} \textsl{ Conformally
%  invariant powers of the {L}aplacian. {I}. {E}xistence},
%\newblock J. London Math. Soc. (2) \textbf{ 46}(3), 557--565 (1992).
%%
%%
%\bibitem[GL91]{GL}
%\msn{1112625}{C.R. Graham, and J.M. Lee,} \textsl{Einstein metrics with prescribed conformal infinity on the ball}, Adv. Math. {\bf 87}, 186--225 (1991).
%%
%%
%\bibitem[GR17]{GrahamRiechert}
%C.R. Graham and N. Riechert, \textsl{
%Higher-dimensional Willmore energies via minimal submanifold asymptotics},
%\arxiv{1704.03852}.
%%
%%
%\bibitem[GrW99]{GrahamWitten}
%\msn{1682674}{C.~R. Graham and E.~Witten,} \textsl{ Conformal anomaly of submanifold
%  observables in {A}d{S}/{CFT} correspondence},
%\newblock Nuclear Phys. B \textbf{ 546}(1-2), 52--64 (1999),  \hepth{9901021}.
%
%%
%%\bibitem[GZ03]{GZ}
%%\msn{1965361}{C.R.  Graham and M. Zworski,}
%%\textsl{Scattering matrix in conformal geometry},
%%Invent. Math. {\bf 152} (1), 89--118 (2003),
%%\maths{0109089}.
%
%\bibitem[Gra03]{Grant}
%D.~Grant,
%\newblock
%  \href{http://www.math.auckland.ac.nz/mathwiki/images/5/51/GrantMSc.pdf}{\it A
%  conformally invariant third order Neumann-type operator for hypersurfaces},
%\newblock Master's thesis, University of Auckland, New Zealand, 2003.
%
%
%
%
%%\bibitem[GKP98]{Gubser}
%%\msn{1630766}{S. S.  
%%Gubser, I. R.  Klebanov,
%%A. M. Polyakov},
%%\textsl{Gauge theory correlators from non-critical string theory},
%%Phys. Lett. B {\bf 428} (1998), 105--114, 
%%\hepth{9802109}.
%
%
%%\bibitem[Guv05]{Guven}
%% \msn{2169323}
%%    {J.~Guven,} \textsl{ Conformally invariant bending energy for hypersurfaces},
%%\newblock J. Phys. A \textbf{ 38}(37), 7943--7955 (2005), \href{http://lanl.arxiv.org/abs/cond-mat/0507320}{arXiv:cond-mat/0507320}.
%
%\bibitem[HS98]{Henningson}
%  \msn{1644988}{M.\ Henningson and K.\ Skenderis,}
%  \textsl{The Holographic Weyl anomaly},
%  JHEP {\bf 9807}, 023 (1998),
% \hepth{hep-th9806087}. 
%
%%\bibitem[LL51]{Landau} \msn{0475345 }{L.D. Landau and E.M Lifshitz,} \textsl{The Classical Theory of Fields}, Course of Theoretical Physics Series, Volume 2, 4 Edition, Butterworth-Heinemann,  Oxford, 1980. 
%
%%\bibitem[LMP01]{grtensor} K. Lake, P. Musgrave and D.  Pollney, GRTensorII, http://grtensor.phy.queensu.ca/, 2001, Maple and Mathematica package. 
%%
%%\bibitem[LM13]{Lewkowycz}
%%\msn{3106348}{A. Lewkowycz and J. Maldacena,} \textsl{Generalized gravitational entropy}, JHEP {\bf 1308} (2013), 090, \arxiv{1304.4926}.
%
%\bibitem[LeB82]{LeB-Heaven}
%\msn{0652038}{C. R. LeBrun,} {\em~$\mathscr{H}$-Space with a Cosmological Constant},
%Proc. R. Soc. Lond.  A {\bf  380},
% 171--185  (1982). 
%
%
%
%%\bibitem[LN74]{Loewner}
%%\msn{0358078}{C.~Loewner and L.~Nirenberg,}
%%\newblock Partial differential equations invariant under conformal or
%%  projective transformations,
%%\newblock in \textsl{ Contributions to analysis (a collection of papers
%%  dedicated to {L}ipman {B}ers)}, pages 245--272, Academic Press, New York,
%%  1974.
%
%%\bibitem[Mal07]{Malchiodi}
%% \msn{2366902}{A. Malchiodi},
%%     \textsl{Conformal metrics with constant {$Q$}-curvature},
%%  SIGMA
%%     {\bf 3},.Paper 120, 11 (2007),
%%    \arxiv{0712.2123} 
%    
%
%  
%  \bibitem[Mal98]{Mal} \msn{1633016}{J.\ Maldacena,} \textsl{The large~$N$ limit of superconformal
%    field theories and supergravity},  Adv.\ Theor.\ Math.\ Phys.\  {\bf 2}, 231--252, (1998), \hepth{9711200}. 
%

%
%\bibitem[Maz91]{MazzeoC}
%\msn{1142715}{R.~Mazzeo,} \textsl{ Regularity for the singular {Y}amabe problem},
%\newblock Indiana Univ. Math. J. \textbf{ 40}(4), 1277--1299 (1991).
%
%\bibitem[Ndi08]{Ndi08}
%\msn{247136}{C. B.  Ndiaye},
%\textsl{Conformal metrics with constant $Q$-curvature for manifolds with boundary}, Comm. Anal. Geom. {\bf 16} (2008),  1049--1124.
%
%
%%\bibitem[Ndi09]{Ndi09}\msn{2485478}{C. B.  Ndiaye}, \textsl{Constant T-curvature conformal metrics on 4-manifolds with boundary}, Pacific J. Math. {\bf 240} (2009),  151--184,  \arxiv{0708.0732}.
%
%%\bibitem[Ndi11]{Ndi11}
%%\msn{2836061}{C. B. Ndiaye},
%%\textsl{$Q$-curvature flow on 4-manifolds with boundary} Math. Z. {\bf 269} (2011),  83--114, \arxiv{0708.2029}.
%
%
%
%\bibitem[OF03]{Osher}
%\msn{1939127}{S.~Osher and R.~Fedkiw,}
%\newblock \textsl{ Level set methods and dynamic implicit surfaces}, volume 153
%  of \textsl{ Applied Mathematical Sciences},
%\newblock Springer-Verlag, New York, 2003.
%
%
%
%\bibitem[PR84]{ot} \msn{917488}{R.\ Penrose, W.\ Rindler,} {\em Spinors and
%  space-time. Vol. 1. Two-spinor calculus and relativistic
%  fields}, Cambridge Monographs on Mathematical Physics. Cambridge
%  University Press, Cambridge, 1984. 
%
%
%%\bibitem[PRR15]{Perlmutter}
%%E.~Perlmutter, M.~Rangamani and M.~Rota, \textsl{ {Positivity, negativity, and
%%  entanglement}},
%%\newblock \arxiv{1506.01679}.
%
%
%\bibitem[Pol86]{Polyakov}
%\msn{0834521}{A. M. Polyakov,} \textsl{Fine Structure of Strings}, Nucl. Phys. B {\bf 268} (1986), 406--412.
%
%%%\cite{Shaukat:2009hp}
%%\bibitem[SW10]{Shaukat} 
%%  \msn{2580074}{A.~Shaukat and A.~Waldron,}
%%  \textsl{Weyl's Gauge Invariance: Conformal Geometry, Spinors, Supersymmetry, and Interactions},
%%  Nucl.\ Phys.\ B {\bf 829}, 28--47 (2010), \arxiv{0911.2477}.
%
%\bibitem[Sus95]{Susskind}
%\msn{1355913}{L. Susskind},  \textsl{The World as a Hologram}, J. of Math. Phys. {\bf 36}  6377?6396 (1995), \hepth{9409089}.
%
%
%\bibitem[Riv08]{Riviere}
%\msn{2430975}
%   {T. Rivi{\`e}re},
%    {\em Analysis aspects of {W}illmore surfaces},
%   {Invent. Math.},
%  {\bf 174}, 1--45
%      (2008), \maths{0612526}.
%
%
%
%\bibitem[RT06]{RyuT1} 
%  S.~Ryu and T.~Takayanagi,
%  {\em Holographic derivation of entanglement entropy from AdS/CFT,}
%  Phys.\ Rev.\ Lett.\  {\bf 96}, 181602 (2006),
%  \hepth{0603001}.
%  
%  
%  \bibitem[RT061]{RyuT2} 
%  S.~Ryu and T.~Takayanagi,
%  {\em Aspects of Holographic Entanglement Entropy,}
%  JHEP {\bf 0608}, 045 (2006),
%  \hepth{0605073}.
%
%
%
%
%\bibitem[Sch84]{Schoen}
%\msn{0788292}{R. Schoen},
%\textsl{Conformal deformation of a Riemannian metric to constant scalar curvature}, J. Differential Geom. {\bf 20}, 479--495, 1984.
%
%%\bibitem[She07]{Sheshadri}
%%\msn{2345277}{N. Seshadri}, 
%%\textsl{Volume renormalization for complete Einstein-K\"ahler metrics} Differential Geom. Appl. {\bf 25} (2007),  356--379,
%%\maths{0404455}
%
%
%\bibitem[Sta05]{Stafford}
%R.~Stafford,
%\newblock
%  \href{www.math.auckland.ac.nz/mathwiki/images/c/cf/StaffordMSc.pdf}{\it
%  Tractor Calculus and Invariants for Conformal Sub-Manifolds},
%\newblock Master's thesis, University of Auckland, New Zealand, 2005.
%
%
%
%\bibitem[Tho26]{Thomas} T.Y.\ Thomas, {\em  On conformal geometry},
%  Proc.\ Natl.\ Acad.\ Sci.\ USA
%  {\bf 12}, 352--359 (1926).
%  Cited on pages
%
%
%%\cite{tHooft:1993dmi}
%\bibitem[tHo93]{'t Hooft} 
%  G.~'t Hooft,
%  \textsl{Dimensional reduction in quantum gravity},
%  Conf.\ Proc.\ C {\bf 930308}, 284--296 (1993)
%  [gr-qc/9310026].
%
%\bibitem[Tru68]{Trudinger}
%\msn{0240748}{N.S. Trudinger},
%\textsl{Remarks concerning the conformal deformation of Riemannian structures on compact manifolds}, Ann. Scuola Norm. Sup. Pisa (3) {\bf 22},265--274, 1968.
%
%\bibitem[Vya13]{YuriThesis}
%Y.~Vyatkin,
%\newblock \textsl{
%  \href{http://librarysearch.auckland.ac.nz/UOA2_A:uoa_alma21234030790002091}{Manufacturing
%  conformal invariants of hypersurfaces}},
%\newblock PhD thesis, University of Auckland, 2013.
%
%%\bibitem[Wal84]{Wald}
%%\msn{0757180}{R. M. Wald,} \textsl{General Relativity}, University of Chicago Press, 2010. 
%
%
%\bibitem[Wil65]{Willmore}
%\msn{0202066}{T.J. Willmore},
%    {\em Note on embedded surfaces},
%  An. \c Sti. Univ. ``Al. I. Cuza'' Ia\c si Sec\c t. I a Mat. (N.S.) {\bf 11B}, 493--496 (1965).
%
%
%
%%\bibitem[Wit98]{Witten}
%%\msn{1633012}{E. Witten}, \textsl{Anti de Sitter space and holography},
%%Adv. Theor. Math. Phys. {\bf 2} (1998), 253--291,
%%\hepth{9802150}.
%
%\bibitem[Yam60]{Yamabe}
%\msn{0125546}{H. Yamabe}, 
%\textsl{On a deformation of Riemannian structures on compact manifolds}, Osaka Journal of Mathematics, {\bf 12}, 21--37,
%1960.
%
\end{thebibliography}
\end{document}